\documentclass[10pt]{article}
\usepackage{xcolor}
\usepackage{amssymb}
\usepackage{cite}
\newlength{\minitwocolumn}
\font\teneufm=eufm10
\font\seveneufm=eufm7
\font\fiveeufm=eufm5
\newfam\eufmfam
\textfont\eufmfam=\teneufm
\scriptfont\eufmfam=\seveneufm
\scriptscriptfont\eufmfam=\fiveeufm

\makeatletter
\renewcommand{\@biblabel}[1]{\textsuperscript{#1}}
\makeatother
\renewcommand{\theequation}{\arabic{equation}}
\newtheorem{thm}{Theorem}[section]

\newtheorem{dfn}[thm]{Definition}
\newtheorem{prop}[thm]{Proposition}
\newtheorem{lem}[thm]{Lemma}

\usepackage{amsmath,amssymb}
\makeatletter
\def\@cite#1{\textsuperscript{#1}}
\makeatother
\title{
\large{\bf
Quadratic relations of the deformed $W$-superalgebra 
${\cal W}_{q, t}(\mathfrak{sl}(2|1))$}}
\author{Takeo KOJIMA}
\begin{document}

\maketitle

\begin{center}
{\it Department of Mathematics and Physics, Faculty of Engineering, Yamagata University,\\
Jonan 4-chome 3-16, Yonezawa 992-8510, JAPAN}
\end{center}

\begin{abstract}
We revisit the free field construction of the deformed $W$-superalgebras ${\cal W}_{q, t}(\mathfrak{sl}(2|1))$ 
by J.\ Ding and B.\ Feigin, {\it Contemp.\ Math.\ }{\bf 248}, 83-108 (1998), where the basic $W$-current and screening currents have been found.  
In this paper we introduce higher $W$-currents and obtain a closed set of quadratic relations among them. 
These relations are independent of the choice of 
Dynkin diagrams for the superalgebra $\mathfrak{sl}(2|1)$,
though the screening currents are not. This allows us to define ${\cal W}_{q, t}(\mathfrak{sl}(2|1))$ by generators and relations.
\end{abstract}

\section{Introduction}
\label{Section1}

The deformed $W$-algebra ${\cal W}_{q,t}(\mathfrak{g})$ is a two parameter deformation
of the classical $W$-algebra ${\cal W}(\mathfrak{g})$, including 
the $W$-algebra and the $q$-Poisson $W$-algebra as special cases.
Shiraishi et al.\cite{Shiraishi-Kubo-Awata-Odake} obtained a free field realization of the deformed Virasoro algebra
${\cal W}_{q,t}(\mathfrak{sl}(2))$, which is a one-parameter deformation of the Virasoro algebra,
to construct a deformation of the correspondence between conformal field theory and the Calogero-Sutherland model.
Subsequently this work has been extended in various directions
\cite{Awata-Kubo-Odake-Shiraishi, Feigin-Frenkel, Brazhnikov-Lukyanov, Hara-Jimbo-Konno-Odake-Shiraishi, Frenkel-Reshetikhin, 
Sevostyanov, Odake, Feigin-Jimbo-Mukhin-Vilkoviskiy} concerning the deformed $W$-algebra ${\cal W}_{q,t}({\mathfrak g})$. 
In comparison with the conformal case, the theory of
deformed $W$-algebras is still not fully developed and understood. 
For that matter it is worthwhile to concretely construct ${\cal W}_{q,t}({\mathfrak g})$ in each case.

In a way different 
from the above papers, Ding and Feigin\cite{Ding-Feigin} constructed free field realizations of 
the deformed $W$-algebras ${\cal W}_{q,t}(\mathfrak{sl}(3))$ and ${\cal W}_{q,t}(\mathfrak{sl}(2|1))$.
They started from a $W$-current $T_1(z)$ given as a sum of three vertex operators $\Lambda_i(z)$ (see (\ref{def:T1(z)})), 
and two screening currents $S_j(w)$ each given by a vertex operator (see (\ref{def:screening})). 
The idea of Ref.\cite{Ding-Feigin} is to determine them simultaneously by demanding that 
$T_1(z)$ and  $S_j(w)$ commute up to a total difference.
Using this method Ding and Feigin obtained three kinds of free field realizations. 
The first realization coincided with that of the deformed $W$-algebra ${\cal W}_{q,t}(\mathfrak{sl}(3))$
obtained earlier in Refs.\cite{Awata-Kubo-Odake-Shiraishi, Feigin-Frenkel}
The second and the third realizations gave objects which were new at that time.
In the conformal limit the second and the third screening charges in Ref.\cite{Ding-Feigin}
satisfy the Serre relations of $U_q(\mathfrak{sl}(2|1))$ associated with the Dynkin diagrams
%%%%%%%%%%%%%%%%%%%%%%%%%%%%%%%%%%%
{\unitlength 0.1in%
\begin{picture}(4.0000,1.0000)(5.5700,-4.0600)%
\special{pn 13}%
\special{ar 607 356 50 50 0.0000000 6.2831853}%
\special{pn 13}%
\special{ar 907 356 50 50 0.0000000 6.2831853}% 
\special{pn 13}%
\special{pa 657 356}%
\special{pa 857 356}%
\special{fp}%
\special{pn 13}%
\special{pa 872 323}%
\special{pa 940 395}%
\special{fp}%
\special{pn 13}%
\special{pa 940 320}%
\special{pa 873 394}%
\special{fp}%
\end{picture}}%
%%%%%%%%%%%%%%%%%%%%%%%%%%%%%%%%
~and
%%%%%%%%%%%%%%%%%%%%%%%%%%%%%%%%
{\unitlength 0.1in%
\begin{picture}(4.0000,1.0000)(5.5700,-4.0600)%
\special{pn 13}%
\special{ar 607 356 50 50 0.0000000 6.2831853}%
\special{pn 13}%
\special{ar 907 356 50 50 0.0000000 6.2831853}%
\special{pn 13}%
\special{pa 657 356}%
\special{pa 857 356}%
\special{fp}%
\special{pn 13}%
\special{pa 642 324}%
\special{pa 573 395}%
\special{fp}%
\special{pa 573 393}%
\special{pa 573 393}%
\special{fp}%
\special{pn 13}%
\special{pa 571 323}%
\special{pa 571 323}%
\special{fp}%
\special{pa 571 325}%
\special{pa 642 393}%
\special{fp}%
\special{pn 13}%
\special{pa 872 323}%
\special{pa 940 395}%
\special{fp}%
\special{pn 13}%
\special{pa 940 320}%
\special{pa 873 394}%
\special{fp}%
\end{picture}}%
%%%%%%%%%%%%%%%%%%%%%%%%%%%%%%%%%%%%%%%%
, respectively.
For that reason they called the algebras generated by the second and the third realizations
the deformed $W$-superalgebra ${\cal W}_{q,t}(\mathfrak{sl}(2|1))$. 
Although the term ``realization'' is used, 
we stress that an abstract definition of ${\cal W}_{q,t}(\mathfrak{sl}(2|1))$ was missing at this stage. 

In this paper we continue the study started in Ding and Feigin\cite{Ding-Feigin}.
We begin by a review of their work. 
Since the derivation of the results given in Ref.\cite{Ding-Feigin} is somewhat sketchy,
we give here the proofs in full detail.
We then introduce higher $W$-currents $T_i(z)$ $(i=1,2,3, \cdots)$,  
and present a closed set of quadratic relations for them. 
In the case of ${\cal W}_{q,t}(\mathfrak{sl}(3))$ 
it is known\cite{Awata-Kubo-Odake-Shiraishi, Feigin-Frenkel} 
that these currents truncate, i.e., $T_3(z)=1$ and $T_i(z)=0$ $(i \geq 4)$, so that the quadratic relations close 
between $T_1(z)$ and $T_2(z)$ (see (\ref{CaseI:quadratic})).
In contrast, such truncation does not take place for the deformed $W$-superalgebra 
${\cal W}_{q,t}(\mathfrak{sl}(2|1))$. 
We show that an infinite number of quadratic relations are satisfied by an infinite number of $T_i(z)$'s
(see (\ref{CaseII:Ti(z)}) and (\ref{CaseII-III:quadratic})).
Moreover, these quadratic relations do not depend on the choice of the Dynkin diagram for the superalgebra $\mathfrak{sl}(2|1)$, 
even though the screening currents do. 
This leads us to define the algebra ${\cal W}_{q,t}(\mathfrak{sl}(2|1))$ abstractly 
by generators  $T_i(z)$ ($i=1,2,3, \cdots$) and defining relations.
This is the main result of this paper.

The text is organized as follows.
In Section \ref{Section2},  
we introduce our notation
and review Ding-Feigin's construction of
the deformed $W$-algebra ${\cal W}_{q,t}(\mathfrak{sl}(3))$ and ${\cal W}_{q,t}(\mathfrak{sl}(2|1))$. 
In Section \ref{Section3}, we introduce higher $W$-currents $T_i(z)$ 
and present a closed set of quadratic relations among them. 
We also obtain the $q$-Poisson algebra in the classical limit. 
Section \ref{Section4} is devoted to conclusion and discussion.

\section{Preliminaries}
\label{Section2}

In this section we prepare the notation and review Ding-Feigin construction of
the deformed $W$-algebra\cite{Ding-Feigin}.
Throughout this paper, we fix a real number $r>1$ and a complex number $x$ with $0<|x|<1$.

\subsection{Notation}

In this section we use complex numbers 
$a$, $w$ $(w\neq 0)$, $q$ ($q \neq 0,\pm1$), 
and $p$ with $|p|<1$.
For any integer $n$, define $q$-integer
\begin{eqnarray}
[n]_q=\frac{q^n-q^{-n}}{q-q^{-1}}.
\nonumber
\end{eqnarray}
We use symbols for infinite products
\begin{eqnarray}
(a;p)_\infty=\prod_{k=0}^\infty (1-a p^k),~~~
(a_1,a_2, \ldots, a_N; p)_\infty=
\prod_{i=1}^N (a_i; p)_\infty.
\nonumber
\end{eqnarray}
for complex numbers $a_1, a_2, \ldots, a_N$.
The following standard formulae are useful.
\begin{eqnarray}
\exp\left(-\sum_{m=1}^\infty \frac{1}{m}a^m \right)=1-a,~~~\exp\left(-\sum_{m=1}^\infty \frac{1}{m}\frac{a^m}{1-p^m}\right)=(a;p)_\infty.\nonumber
%\label{eqn:infinite-product}
\end{eqnarray}
We use the elliptic theta function
\begin{eqnarray}
\Theta_p(w)=(p, w, p w^{-1};p)_\infty,~~~
\Theta_p (w_1, w_2, \ldots, w_N)=\prod_{i=1}^N \Theta_p(w_i)
\nonumber
\end{eqnarray}
for complex numbers $w_1, w_2, \ldots, w_N \neq 0$.
Define $\delta(z)$ by the formal series
\begin{eqnarray}
\delta(z)=\sum_{m \in {\mathbf Z}}z^m.
\nonumber%\label{def:delta}
\end{eqnarray}

\subsection{Ding-Feigin realization}

In the following we review the results of Ding and Feigin\cite{Ding-Feigin}. 
We state and rederive their results, giving the full detail and 
correcting some minor mistakes in their paper. 

We introduce the Heisenberg algebra with generators
$a_i(m)$, $Q_i$ $(m \in {\bf Z}, 1\leq i \leq 2)$ satisfying
\begin{eqnarray}
&&
[a_i(m), a_j(n)]=\frac{1}{m}A_{i,j}(m)\delta_{m+n,0}~~
(m, n \neq 0, 1\leq i,j \leq 2),
\nonumber\\
&&
[a_i(0),Q_j]=A_{i,j}(0)~~(1\leq i,j \leq 2).
\nonumber
\end{eqnarray}
The remaining commutators vanish.
We impose the following conditions on the parameters $A_{i,j}(m)\in\mathbf{C}$:
\begin{eqnarray}
&&
A_{1,1}(m)=A_{2,2}(m)=1~(m \neq 0),~~~A_{1,2}(m)=A_{2,1}(-m)~(m \in \mathbf{Z}),
\nonumber\\
&&
A_{1,2}(m)A_{2,1}(m)\neq A_{1,1}(m)A_{2,2}(m)~~(m \in {\mathbf Z}).
\nonumber
\end{eqnarray}
We use the normal ordering symbol $:~:$ that satisfies
\begin{eqnarray}
&&
:a_i(m)a_j(n):=\left\{\begin{array}{cc}
a_i(m)a_j(n)&(m<0),\\
a_j(n)a_i(m)&(m \geq 0)
\end{array}
\right.~~~(m,n \in {\mathbf Z}, 1\leq i,j \leq 2),
\nonumber
\\
&&
:a_i(0)Q_j:=:Q_j a_i(0):=Q_j a_i(0)~~~(1 \leq i,j \leq 2).
\nonumber
\end{eqnarray}

Next, we work on Fock space of the Heisenberg algebra.
Let $T_1(z)$ be a sum of three vertex operators
\begin{align}
&T_1(z)=g_1 \Lambda_1(z)+g_2 \Lambda_2(z)+g_3 \Lambda_3(z),
\label{def:T1(z)}\\
&{\Lambda}_i(z)=e^{\sum_{j=1}^2 \lambda_{i,j}(0)a_j(0)}
:\exp\left(\sum_{j=1}^2 \sum_{m \neq 0}\lambda_{i,j}(m)a_j(m)z^{-m}\right):~~(1\leq i \leq 3).
\nonumber
\end{align}
We call  $T_1(z)$ the basic $W$-current.
We introduce the screening currents 
$S_j(w)$ $(1 \leq j \leq 2)$ as
\begin{eqnarray}
S_j(w)=w^{\frac{1}{2}A_{j,j}(0)}e^{Q_j}w^{a_j(0)}
:\exp\left(\sum_{m \neq 0}s_j(m)a_j(m)w^{-m}\right):
~~(1\leq j \leq 2).
\label{def:screening}
\end{eqnarray}
The complex parameters
$A_{ij}(m)$, $\lambda_{ij}(m), s_j(m)$ and $g_i$ are to be determined through the construction given below.

Quite generally, given two vertex  operators $V(z)$, $W(w)$, 
their product has the form 
\begin{eqnarray}
V(z)W(w)=\varphi_{V,W}\left(z,w\right):V(z)W(w):~~~(|z|\gg |w|),
\nonumber
\end{eqnarray}
with some formal power series $\varphi_{V,W}(z,w)\in\mathbf{C}[[w/z]]$. 
The vertex operators $V(z)$ and $W(w)$ are said to be mutually local if the following two conditions hold.
\begin{align}
&\mathrm{(i)}~~\varphi_{V,W}(z,w)~{\rm and}~\varphi_{W,V}(w,z)~{\rm converge~ to~rational~functions},~~~\nonumber\\
&\mathrm{(ii)}~~\varphi_{V,W}(z,w)=\varphi_{W,V}(w,z).\nonumber
\end{align}

Under this setting, we are going to determine
the $W$-currents $T_1(z)$ 
and the screening currents $S_j(w)$ that
satisfy the following
mutual locality (\ref{assume:mutual-locality}), commutativity (\ref{assume:q-difference}),
and symmetry (\ref{assume:varphiSS}).
\\
\\
{\bf Mutual Locality}~~~$\Lambda_i(z)$ $(1\leq i \leq 3)$ and 
$S_j(w)$ $(1\leq j \leq 2)$ are mutually local, the operator product expansions of their products
have at most one pole and one zero, and
\begin{eqnarray}
\varphi_{\Lambda_i, S_j}(z,w)=\varphi_{S_j, \Lambda_i}(w,z)
={\displaystyle \frac{w-\frac{z}{p_{i,j}}}{w-\frac{z}{q_{i,j}}}}~~~(1\leq i \leq 3, 1\leq j \leq 2).
\label{assume:mutual-locality}
\end{eqnarray}
We allow the possibility $p_{i,j}=q_{i,j}$,  in which case 
$\Lambda_i(z)S_j(w)=S_j(w)\Lambda_i(z)=:\Lambda_i(z)S_j(w):$.

~\\
{\bf Commutativity}~~~$T_1(z)$ commutes with $S_j(w)$ $(1\leq j \leq 2)$
up to a total difference
\begin{eqnarray}
[T_1(z),S_j(w)]=B_j(z)
\left(\delta\left(\frac{q_{j,j}w}{z}\right)-\delta\left(\frac{q_{j+1,j}w}{z}\right)\right)~~~(1\leq j \leq 2)\,,
\label{assume:q-difference}
\end{eqnarray}
with some currents $B_j(z)$ $(1\leq j \leq 2)$.

~\\
{\bf Symmetry}~~~For $\widetilde{S}_j(w)=e^{-Q_j}S_j(w)$ $(1\leq j \leq 2)$, we impose
\begin{eqnarray}
\varphi_{\widetilde{S}_1, \widetilde{S}_2}(w,z)=
\varphi_{\widetilde{S}_2, \widetilde{S}_1}(w,z).
\label{assume:varphiSS}
\end{eqnarray}
For simplicity, we impose further the following conditions.
\begin{align}
&q_{i,j}~~(1\leq i \leq 3, 1\leq j \leq 2)~~{\rm are~distinct}.
\label{assume:p,q0}\\
&\left|\frac{q_{j+1,j}}{q_{j,j}}\right| \neq 1~~~(1\leq j \leq 2),~~~-1<A_{1,2}(0)=A_{2,1}(0)<0.
\label{assume:p,q1}
\end{align}
We introduce the two parameters $x$ and $r$ defined as
\begin{eqnarray}
x^{2r}=\frac{q_{2,1}}{q_{1,1}},~~~r=\left\{\begin{array}{cc}
{{\displaystyle \frac{1}{A_{1,2}(0)+1}}}&(A_{1,1}(0)\neq 1),\\
{\displaystyle -\frac{1}{A_{1,2}(0)}}&(A_{1,1}(0)=1).
\end{array}
\right.
\label{def:2parameters}
\end{eqnarray}
From (\ref{assume:p,q1}) and $q_{i,j}\neq 0$,
we obtain $|x|\neq 0, 1$ and $r>1$.
In this paper, we focus our attention to
\begin{eqnarray}
0<|x|<1,~~~ r>1.\label{regime:x,r}
\nonumber
\end{eqnarray}
For the case of $|x|>1$, we obtain the same results under the change $x \mapsto x^{-1}$.

Consider the following transformations which map operators of form (\ref{def:T1(z)}), (\ref{def:screening}) into operators of the same form.
\\
{(i)}~Rearranging indices
\begin{eqnarray}
\Lambda_{i}(z) \mapsto \Lambda_{i'}(z),~~~S_j(w) \mapsto S_{j'}(w),
\label{def:arrange-suffix}
\end{eqnarray}
where $i \to i'$ is a permutation of the set $1,2,3$
and $j \to j'$ is a permutation of the set $1,2$.\\
{(ii)}~Scaling variables: $\Lambda_i(z) \mapsto \Lambda_i(sz)$ 
$(s \neq 0)$, i.e.
\begin{eqnarray}
\lambda_{i,j}(m) \mapsto s^m \lambda_{i,j}(m)~~~
q_{i,j} \mapsto s q_{i,j},~~~p_{i,j} \mapsto s p_{i,j}~~~(m\neq 0, 1\leq i \leq 3, 1\leq j \leq 2).
\label{def:shift}
\end{eqnarray}
{(iii)}~Scaling free fields :
\begin{eqnarray}
\begin{array}{cc}
\begin{array}{c}
a_j(m) \mapsto \alpha_j(m)^{-1} a_j(m),~~s_j(m) \mapsto \alpha_j(m) s_j(m),\\
\lambda_{i,j}(m) \mapsto \lambda_{i,j}(m)\alpha_j(m)
\end{array}&(m \neq 0, 1\leq i \leq 3, 1 \leq j \leq 2),
\\
A_{i,j}(m) \mapsto \alpha_i(m)^{-1} A_{i,j}(m) \alpha_j(m)
&(m \neq 0, 1\leq i, j \leq 2),
\end{array}
\label{def:scaling-boson}
\end{eqnarray}
where
$\alpha_j(m)\neq 0$ $(1\leq j \leq 2)$ and $\alpha_j(-m)=\alpha_j(m)^{-1}$ $(m>0, 1\leq j \leq 2)$.

In Ding-Feigin's construction\cite{Ding-Feigin}, there are four cases to be considered separately.
\begin{eqnarray}
\begin{array}{cc}
{\rm Case~1}:~A_{1,1}(0)\neq 1,~ A_{2,2}(0)\neq 1,& 
{\rm Case~2}:~A_{1,1}(0)\neq 1,~ A_{2,2}(0)=1,\\
{\rm Case~3}:~A_{1,1}(0)=1,~ A_{2,2}(0)=1,&
{\rm Case~4}:~A_{1,1}(0)=1,~ A_{2,2}(0)\neq 1.
\end{array}
\label{def:4cases}
\end{eqnarray}
From (\ref{def:arrange-suffix}),
either Case 2 or Case 4 can be omitted.
In what follows, we omit Case 4 and
fix the index of $\Lambda_i(z)$ and $S_j(w)$.

\begin{thm}\label{thm:II-1}
Assume that conditions
(\ref{assume:mutual-locality}), (\ref{assume:q-difference}), 
(\ref{assume:varphiSS}),  (\ref{assume:p,q0}) and (\ref{assume:p,q1})
hold.
Then, up to transformations
\eqref{def:arrange-suffix}, \eqref{def:shift} and \eqref{def:scaling-boson}, 
the parameters  $p_{i,j}, q_{i,j}$, $A_{i,j}(m)$, $s_i(m)$, $\lambda_{i,j}(m)$, $g_i$, and the current $B_j(m)$
are uniquely determined as follows.
\\
$\bullet$~In all Cases 1, 2, and 3, 
\begin{eqnarray}
&&
p_{1,2}=q_{1,2},~~~p_{3,1}=q_{3,1},~~~s_j(m)=1~~~(m>0, 1\leq j \leq 2),
\label{CaseI-II-III:p,q,s(m)}
\\
&&
A_{1,1}(m)=A_{2,2}(m)=1,~~~A_{2,1}(m)=A_{1,2}(-m)~~~(m \neq 0),
\nonumber%\label{CaseI-II-III:A(m)}\\
\\
&&
B_j(z)=g_j\left(\frac{q_{j,j}}{p_{j,j}}-1\right):\Lambda_j(z)S_j(q_{j,j}^{-1}z):~~~(1\leq j \leq 2).
\label{CaseI-II-III:B}
\end{eqnarray}
\\
$\bullet$~In Case 1,
\begin{eqnarray}
&&A_{1,1}(0)=A_{2,2}(0)=\frac{2(r-1)}{r},~~~A_{1,2}(0)=A_{2,1}(0)=-\frac{r-1}{r},~~~
g_1=g_2=g_3,
\label{CaseI:A(0),g}
\\
&&
q_{i,i}=x^{i-1},~q_{i+1,i}=x^{2r+i-1},~p_{i,i}=x^{2r+i-3},~p_{i+1,i}=x^{i+1}~~~(1\leq i \leq 2),
\label{CaseI:p,q}\\
&&
A_{1,2}(m)=-\frac{[m]_x}{[2m]_x}~~~(m \neq 0),
~~~s_j(m)=-\frac{[(r-1)m]_x[2m]_x}{
[rm]_x[m]_x}~~~(m<0, 1\leq j \leq 2),
\label{CaseI:A(m)-s(m)}
\end{eqnarray}
\begin{eqnarray}
&&
\lambda_{i,j}(0)=\frac{2r}{3}\log x \times
\left\{
\begin{array}{cc}
j& (1\leq j<i \leq 3), \\
j-3 &(1\leq i \leq j \leq 2),
\end{array}
\right.
\label{CaseI:lambda(0)}
\\
&&
\frac{\lambda_{i,j}(m)}{s_j(m)}=-\frac{[rm]_x}{[3m]_x}(x-x^{-1})\times
\left\{
\begin{array}{cc}
x^{(r+2)m}[jm]_x& (1\leq j<i \leq 3), \\
x^{(r-1)m}[(j-3)m]_x&(1\leq i \leq j \leq 2)
\end{array}
\right.~~(m \neq 0).
\label{CaseI:lambda(m)}
\end{eqnarray}
$\bullet$~In Case 2,
\begin{eqnarray}
&&A_{1,1}(0)=\frac{2(r-1)}{r},~A_{2,2}(0)=1,~~A_{1,2}(0)=A_{2,1}(0)=-\frac{r-1}{r},
~~g_1=g_2,~~g_3=[r-1]_x g_1,
\label{CaseII:A(0),g}
\\
&&
q_{i,i}=x^{i-1},~q_{i+1,i}=x^{2r+i-1},~p_{i,i}=x^{2r+i-3},~p_{i+1,i}=x^{-2r+5+(2r-3)i}~~~(1\leq i \leq 2),
\label{CaseII:p,q}
\\
&&
A_{1,2}(m)=\left\{
\begin{array}{cc}
{\displaystyle
-\frac{[(r-1)m]_x}{[rm]_x}}
&(m>0),\\
{\displaystyle
-\frac{[m]_x}{[2m]_x}}&(m<0),
\end{array}
\right.~
s_1(m)=
{\displaystyle
-\frac{[(r-1)m]_x[2m]_x}{
[rm]_x[m]_x}},~
s_2(m)=
-1~(m<0),
\label{CaseII:A(m)-s(m)}
\end{eqnarray}
\begin{eqnarray}
&&
\lambda_{i,j}(0)=\frac{2r}{r+1}\log x \times
\left\{
\begin{array}{cc}
-r& (i,j)=(1,1),\\
1-r &(i,j)=(1,2), (2,2),\\
1&(i,j)=(2,1), (3,1),\\
2&(i,j)=(3,2),
\end{array}
\right.\label{CaseII:lambda(0)}
\\
&&
\frac{\lambda_{i,j}(m)}{s_j(m)}=\frac{[rm]_x}{[(r+1)m]_x}(x-x^{-1})\times
\left\{
\begin{array}{cc}
x^{(r-1)m}[r m]_x& (i,j)=(1,1),
\\
x^{(r-1)m}[(r-1)m]_x&(i,j)=(1,2), (2,2),
\\
-x^{2rm}[m]_x&(i,j)=(2,1), (3,1),
\\
-x^{2rm}[2m]_x&(i,j)=(3,2)
\end{array}
\right.~~(m \neq 0).
\label{CaseII:lambda(m)}
\end{eqnarray}
$\bullet$~In Case 3,
\begin{eqnarray}
&&A_{1,1}(0)=A_{2,2}(0)=1,~~A_{1,2}(0)=A_{2,1}(0)=-\frac{1}{r},
~~g_2=[r-1]_xg_1,~~~g_3=g_1,
\label{CaseIII:A(0),g}
\\
&&
q_{i,i}=x^{(r-1)(i-1)},~q_{i+1,i}=x^{r+1+(r-1)i},~p_{i,i}=p_{i+1,i}
=x^{3r-5+(-r+3)i}~~~(1\leq i \leq 2),
\label{CaseIII:p,q}
\\
&&
A_{1,2}(m)=-\frac{[m]_x}{[rm]_x}~~(m \neq 0),~~~
s_1(m)=s_2(m)=-1~~(m<0),
\label{CaseIII:A(m)-s(m)}
\\
&&
\lambda_{i,j}(0)=\frac{2r}{r+1}\log x \times
\left\{
\begin{array}{cc}
-r& (i,j)=(1,1),\\
-1&(i,j)=(1,2), (2,2),\\
1&(i,j)=(2,1), (3,1),\\
r&(i,j)=(3,2),
\end{array}
\right.\label{CaseIII:lambda(0)}
\\
&&
\frac{\lambda_{i,j}(m)}{s_j(m)}=\frac{[rm]_x}{[(r+1)m]_x}(x-x^{-1})\times
\left\{
\begin{array}{cc}
x^{(r-1)m}[r m]_x& (i,j)=(1,1),\\
x^{(r-1)m}[m]_x&(i,j)=(1,2), (2,2),\\
-x^{2rm}[m]_x&(i,j)=(2,1), (3,1),\\
-x^{2rm}[r m]_x&(i,j)=(3,2)
\end{array}
\right.~~(m \neq 0).
\label{CaseIII:lambda(m)}
\end{eqnarray}
Conversely, if the parameters
are chosen as above then
(\ref{assume:mutual-locality}), (\ref{assume:q-difference}), and
(\ref{assume:varphiSS}) are satisfied.
\end{thm}

\begin{prop}
\label{prop:II-2}
The $\Lambda_i(z)$'s satisfy the commutation relations
\begin{eqnarray}
\Lambda_k(z_1)\Lambda_l(z_2)=-\frac{z_2}{z_1}
\frac{\Theta_{x^{2s}}\left(x^2\frac{z_2}{z_1},~~x^{2s-2r}\frac{z_2}{z_1},~~x^{2s+2r-2}\frac{z_2}{z_1}\right)}{
\Theta_{x^{2s}}\left(x^2\frac{z_1}{z_2},~~x^{2s-2r}\frac{z_1}{z_2},~~x^{2s+2r-2}\frac{z_1}{z_2}\right)
}\Lambda_l(z_2)\Lambda_k(z_1)~~~(1\leq k,l \leq 3)\,,
\label{CaseI-II-III:vertex}
\end{eqnarray}
where 
\begin{align*}
s=\begin{cases}
3 & \text{for Case 1},\\
r+1 & \text{for Cases 2 and 3}.\\
\end{cases}
\end{align*}
We understand (\ref{CaseI-II-III:vertex}) in the sense of analytic continuation.
\end{prop}

\begin{prop}
\label{prop:II-3}
~For Case 1, $S_j(w)$ satisfy
\begin{eqnarray}
&&
S_1(w_1)S_2(w_2)=\left(\frac{w_1}{w_2}\right)^{1+\frac{1}{r}} 
\frac{\Theta_{x^{2r}}
\left(x^{-1} \frac{w_2}{w_1}\right)}{\Theta_{x^{2r}}\left(x^{-1}\frac{w_1}{w_2}\right)}
S_2(w_2)S_1(w_1),
\label{CaseI:S1S2}
\\
&&
S_j(w_1)S_j(w_2)=-\left(\frac{w_1}{w_2}\right)^{1-\frac{2}{r}}
\frac{\Theta_{x^{2r}}\left(x^{2}\frac{w_2}{w_1}\right)}{
\Theta_{x^{2r}}\left(x^{2}\frac{w_1}{w_2}\right)}S_j(w_2)S_j(w_1)~~~(1\leq j \leq 2).
\label{CaseI:SjSj}
\end{eqnarray}

For Case 2, $S_j(w)$ satisfy (\ref{CaseI:S1S2}), (\ref{CaseI:SjSj}) for $j=1$, and
\begin{eqnarray}
S_2(w_1)S_2(w_2)=-S_2(w_2)S_2(w_1).
\label{CaseII:SjSj}
\end{eqnarray}

For Case 3, $S_j(w)$ satisfy 
\begin{eqnarray}
&&
S_1(w_1)S_2(w_2)=\left(\frac{w_1}{w_2}\right)^{-\frac{1}{r}} 
\frac{\Theta_{x^{2r}}
\left(x^{r+1} \frac{w_2}{w_1}\right)}{\Theta_{x^{2r}}\left(x^{r+1}\frac{w_1}{w_2}\right)}
S_2(w_2)S_1(w_1),
\label{CaseIII:S1S2}
\\
&&
S_j(w_1)S_j(w_2)=-S_j(w_2)S_j(w_1)~~(1\leq j \leq 2).
\label{CaseIII:SjSj}
\end{eqnarray}

We understand (\ref{CaseI:S1S2}), (\ref{CaseI:SjSj}), (\ref{CaseII:SjSj}), (\ref{CaseIII:S1S2}), and (\ref{CaseIII:SjSj}) on the analytic continuation.
\end{prop}
In fact, the stronger relation
\begin{eqnarray}
S_j(w_1)S_j(w_2)=(w_1-w_2):S_j(w_1)S_j(w_2):
\nonumber
\end{eqnarray}
holds instead of (\ref{CaseII:SjSj}) and (\ref{CaseIII:SjSj}).
This means that the screening currents 
$S_2(w)$ in Case 2 and $S_1(w), S_2(w)$ in Case 3
are ordinary fermions.

We give a few words about Ref.\cite{Ding-Feigin}.
Apart from some typos, the essential content is correct.
We believe that their assumption
$A_{1,2}(m)=A_{2,1}(m)$ is a misprint.
The free field realizations of $\Lambda_i(z)$ $(1\leq i \leq 3)$ 
were not completely constructed.
In Case 1, all $\Lambda_i(z)$ $(1\leq i \leq 3)$ were not constructed, in Case 2, $\Lambda_2(z)$ was not constructed,
and in Cass 3, $\Lambda_i(z)$ $(i=1,2)$ were not
constructed.
The exchange relations
(\ref{CaseI-II-III:vertex}) 
were only partially investigated.
In Case 1, all exchange relations were not investigated,
in Case 2, only $k=l=1$ and $k=l=3$,
and in Case 3, only $k=l=3$ was investigated.

\subsection{Proof of Theorem \ref{thm:II-1}}

In this section we show Theorem \ref{thm:II-1} and
Proposition \ref{prop:II-3}.
\begin{lem}~For $\Lambda_i(z)$ and $S_j(w)$, we obtain
\begin{eqnarray}
&&
\varphi_{\Lambda_i, S_j}(z, w)
=e^{\sum_{k=1}^2 \lambda_{i,k}(0)A_{k,j}(0)}
\exp\left(\sum_{k=1}^2 \sum_{m=1}^\infty
\frac{1}{m}\lambda_{i,k}(m)A_{k,j}(m)s_j(-m)\left(\frac{w}{z}\right)^m
\right),
\label{normal1}
\\
&&
\varphi_{S_j, \Lambda_i}(w, z)
=\exp\left(\sum_{k=1}^2 \sum_{m=1}^\infty
\frac{1}{m}s_j(m)A_{j,k}(m)\lambda_{i,k}(-m)\left(\frac{z}{w}\right)^m
\right)~~~(1\leq i \leq 3, 1\leq j \leq 2),
\label{normal2}
\\
&&
\varphi_{\widetilde{S}_k,\widetilde{S}_l}(w_1, w_2)=
\exp\left(\sum_{m=1}^\infty \frac{1}{m}s_k(m)A_{k,l}(m)s_l(-m)\left(\frac{w_2}{w_1}\right)^m \right)~~~(1\leq k, l \leq 2),
\label{normal3}
\end{eqnarray}
\begin{eqnarray}
\varphi_{\Lambda_k, \Lambda_l}(z_1, z_2)=
\exp\left(\sum_{i, j=1}^2 \sum_{m=1}^\infty
\frac{1}{m}\lambda_{k,i}(m)A_{i,j}(m)\lambda_{l,j}(-m)\left(\frac{z_2}{z_1}\right)^m \right)~~~(1\leq k, l \leq 3).
\label{normal4}
\end{eqnarray}
\end{lem}
{\it Proof.}~
Using the standard formula
\begin{eqnarray}
e^Ae^B=e^{[A,B]}e^Be^A~~~([[A,B],A]=0~{\rm and}~[[A,B],B]=0),
\nonumber
\end{eqnarray}
we obtain (\ref{normal1}), (\ref{normal2}), (\ref{normal3}), and (\ref{normal4}).~~
$\Box$

\begin{lem}
\label{lem:II-5}~
Mutual locality (\ref{assume:mutual-locality}) holds
if and only if (\ref{eqn:lambda(0)})
and (\ref{eqn:lambda(m)}) are satisfied
\begin{eqnarray}
&&
\sum_{k=1}^2\lambda_{i,k}(0)A_{k,j}(0)=\log\left(\frac{q_{i,j}}{p_{i,j}}\right)
~~~(1\leq i \leq 3, 1\leq j \leq 2),
\label{eqn:lambda(0)}
\\
&&
\sum_{k=1}^2
\lambda_{i,k}(m)A_{k,j}(m)s_j(-m)=q_{i,j}^m-p_{i,j}^m
~~~(m \neq 0, 1\leq i \leq 3, 1\leq j \leq 2).
\label{eqn:lambda(m)}
\end{eqnarray}
\end{lem}
{\it Proof.}~~
Considering (\ref{normal1}), (\ref{normal2}), and the expansions
\begin{eqnarray}
&&
{\displaystyle \frac{w-p_{i,j}^{-1} z}{w-q_{i,j}^{-1} z}}=
\exp\left(\log\left(\frac{q_{i,j}}{p_{i,j}}\right)-\sum_{m=1}^\infty \frac{1}{m}(p_{i,j}^m-q_{i,j}^m)\left(\frac{w}{z}\right)^m\right)
~~~(|z| \gg |w|),
\label{expansion1}\\
&&
\frac{w-p_{i,j}^{-1} z}{w-q_{i,j}^{-1} z}=
\exp\left(-\sum_{m=1}^\infty \frac{1}{m}(p_{i,j}^{-m}-q_{i,j}^{-m})\left(\frac{z}{w}\right)^{m}\right)
~~~(|w| \gg |z|),
\label{expansion2}
\end{eqnarray}
we obtain
(\ref{eqn:lambda(0)}) and (\ref{eqn:lambda(m)})
from (\ref{assume:mutual-locality}).

Conversely, if we assume (\ref{eqn:lambda(0)}) and 
(\ref{eqn:lambda(m)}), we obtain (\ref{assume:mutual-locality})
from (\ref{normal1}), (\ref{normal2}), (\ref{expansion1}), and (\ref{expansion2}).~~
$\Box$

From linear equations (\ref{eqn:lambda(0)}) and (\ref{eqn:lambda(m)}),
$\lambda_{i,j}(m)$ are expressed in terms of the other parameters.

\begin{lem}\label{lem:II-6}~
We assume (\ref{assume:mutual-locality}) and (\ref{assume:p,q0}).
The commutativity (\ref{assume:q-difference}) holds
if and only if $p_{1,2}=q_{1,2}$, $p_{3,1}=q_{3,1}$,
(\ref{assume:start1}), (\ref{assume:start2}), and (\ref{def:current-B}) are satisfied, where
\begin{eqnarray}
&&
q_{j,j}^{\frac{1}{2}A_{j,j}(0)}:\Lambda_j(z)S_j
\left(q_{j,j}^{-1}z\right):=
q_{j+1,j}^{\frac{1}{2}A_{j,j}(0)}:\Lambda_{j+1}(z)
S_j\left(q_{j+1,j}^{-1} z \right):
~~~(1\leq j \leq 2),
\label{assume:start1}
\\
&&
\frac{g_{j+1}}{g_j}=-\left(\frac{q_{j+1,j}}{q_{j,j}}\right)^{\frac{1}{2}A_{j,j}(0)}
\frac{\frac{q_{j,j}}{p_{j,j}}-1}{
\frac{~q_{j+1,j}}{p_{j+1,j}}-1~}~~~(1\leq j \leq 2),
\label{assume:start2}
\\
&&
B_j(z)=g_j\left(\frac{q_{j,j}}{p_{j,j}}-1\right):\Lambda_j(z)S_j(q_{j,j}^{-1}z):~~~(1\leq j \leq 2).
\label{def:current-B}
\end{eqnarray}
\end{lem}
{\it Proof.}~~
From (\ref{assume:mutual-locality}), we obtain
\begin{eqnarray}
[\Lambda_i(z),S_j(w)]=\left(\frac{q_{i,j}}{p_{i,j}}-1\right)
\delta\left(\frac{q_{i,j}w}{z}\right):\Lambda_i(z)S_j(q_{i,j}^{-1}z):
~~(1\leq i \leq 3, 1\leq j \leq 2).\label{eqn:lambda-S}
\end{eqnarray}

Considering (\ref{assume:p,q0}) and (\ref{eqn:lambda-S}), we know that
(\ref{assume:q-difference}) holds 
if and only if $p_{1,2}=q_{1,2}$, $p_{3,1}=q_{3,1}$, and
\begin{eqnarray}
B_j(z)=g_j\left(\frac{q_{j,j}}{p_{j,j}}-1\right):\Lambda_j(z)S_j(q_{j,j}^{-1}z):
=-g_{j+1}\left(\frac{q_{j+1,j}}{p_{j+1,j}}-1\right):\Lambda_{j+1}(z)S_j(q_{j+1,j}^{-1}z):~(1\leq j \leq 2)
\label{eqn:current-B}
\end{eqnarray}
are satisfied.
(\ref{eqn:current-B}) holds
if and only if (\ref{assume:start1}), (\ref{assume:start2}), and (\ref{def:current-B}) are satisfied.
Hence, we obtain this lemma.~~~
$\Box$

We use the abbreviation $h_{k,l}(w)$ $(1\leq k,l \leq 2)$,
\begin{eqnarray}
h_{k,l}\left(\frac{w_2}{w_1}\right)=
\varphi_{\widetilde{S}_k, \widetilde{S}_l}(w_1, w_2).
\label{def:h}
\end{eqnarray}

\begin{lem}\label{lem:II-7}~
We assume (\ref{assume:mutual-locality}) and (\ref{assume:start1}).
Then, $h_{k,l}(w)$ in (\ref{def:h}) satisfy the $q$-difference equations
\begin{eqnarray}
&&
\begin{array}{cc}
{\displaystyle \frac{h_{1,2}(q_{1,1}w)}{h_{1,2}(q_{2,1}w)}
=\frac{q_{2,2}}{p_{2,2}}\left(\frac{q_{1,1}}{q_{2,1}}\right)^{A_{1,2}(0)}
\frac{1-p_{2,2}w}{1-q_{2,2}w}},
&
{\displaystyle \frac{h_{1,2}\left(q_{2,2}^{-1}w \right)}{h_{1,2}\left(q_{3,2}^{-1} w \right)}
=\frac{1-q_{2,1}^{-1} w}{1-p_{2,1}^{-1}w}},
\\
{\displaystyle
\frac{h_{2,1}(q_{3,2}w)}{h_{2,1}(q_{2,2}w)}
=\frac{q_{2,1}}{p_{2,1}}\left(\frac{q_{3,2}}{q_{2,2}}\right)^{A_{1,2}(0)}
\frac{1-p_{2,1}w}{1-q_{2,1}w}},&
{\displaystyle
\frac{h_{2,1}\left(q_{1,1}^{-1} w \right)}{h_{2,1}\left(q_{2,1}^{-1} w \right)}
=\frac{1-p_{2,2}^{-1}w}{1-q_{2,2}^{-1} w}},
\end{array}
\label{eqn:q-diff1}
\end{eqnarray}
\begin{eqnarray}
&&
\begin{array}{c}
{\displaystyle
\frac{1-(p_{j,j}w)^{-1}}{1-(q_{j,j}w)^{-1}}
h_{j,j}\left(q_{j,j}^{-1}w\right)=
\frac{1-(p_{j+1,j}w)^{-1}}{1-(q_{j+1, j}w)^{-1}}
h_{j,j}\left(q_{j+1,j}^{-1}w \right)},
\\
{\displaystyle
\left(\frac{q_{j+1,j}}{q_{j,j}}\right)^{A_{j,j}(0)-1}
\frac{p_{j+1,j}}{p_{j,j}}
\frac{1-p_{j,j}w}{1-q_{j,j}w}h_{j,j}\left(q_{j,j}w\right)=
\frac{1-p_{j+1,j}w}{1-q_{j+1,j}w}h_{j,j}\left(q_{j+1,j}w\right)}
\end{array}
~(1\leq j \leq 2).
\label{eqn:hii-1}
\end{eqnarray}
\end{lem}
{\it Proof.}~
Multiplying (\ref{assume:start1}) by the screening currents on the left or right and considering the normal orderings, we obtain
(\ref{eqn:q-diff1}) and (\ref{eqn:hii-1})
as necessary conditions.~~~
$\Box$

\begin{lem}\label{lem:II-8}~
The relation (\ref{eqn:p,q1}) holds if (\ref{assume:mutual-locality}), (\ref{assume:varphiSS}), (\ref{assume:p,q1}), and (\ref{assume:start1}) are satisfied, where
\begin{eqnarray}
&&
q_{1,1}=s,~~q_{2,2}=sx^{(1+A_{1,2}(0))r},
~~q_{2,1}=sx^{2r},~~q_{3,2}=sx^{(3+A_{1,2}(0))r},\nonumber\\
&&p_{2,1}=sx^{2(1+A_{1,2}(0))r},~~p_{2,2}=sx^{(1-A_{1,2}(0))r},\label{eqn:p,q1}
\\
&&
s_1(m)A_{1,2}(m)s_2(-m)=s_2(m)A_{2,1}(m)s_1(-m)
=-\frac{[A_{1,2}(0) r m]_x}{[rm]_x}~~(m>0).
\nonumber
\end{eqnarray}
\end{lem}
{\it Proof.}~
From lemma \ref{lem:II-7}, we obtain
(\ref{eqn:q-diff1}).
From (\ref{normal3}) and (\ref{def:h}), the constant term of $h_{k,l}(w)$ is 1.
Comparing the Taylor expansions for both sides of 
(\ref{eqn:q-diff1}), we obtain
\begin{eqnarray}
\frac{q_{2,2}}{p_{2,2}}\left(\frac{q_{1,1}}{q_{2,1}}\right)^{A_{1,2}(0)}=1,~~
\frac{q_{2,1}}{p_{2,1}}\left(\frac{q_{3,2}}{q_{2,2}}\right)^{A_{1,2}(0)}
=1.\label{eqn:p,q2}
\end{eqnarray}
Upon the specialization (\ref{eqn:p,q2}), 
we obtain solutions of
(\ref{eqn:q-diff1}) as
\begin{eqnarray}
h_{1,2}(w)=\exp\left(-\sum_{m=1}^\infty
\frac{1}{m}\frac{\left(\frac{p_{2,2}}{q_{1,1}}\right)^m-
\left(\frac{q_{2,2}}{q_{1,1}}\right)^m}{1-\left(\frac{q_{2,1}}{q_{1,1}}\right)^m}w^m\right)
=
\exp\left(-\sum_{m=1}^\infty
\frac{1}{m}\frac{\left(\frac{q_{3,2}}{p_{2,1}}\right)^m-
\left(\frac{q_{3,2}}{q_{2,1}}\right)^m}{1-\left(\frac{q_{3,2}}{q_{2,2}}\right)^m}w^m\right),
\label{eqn:h12-1}
\\
h_{2,1}(w)=
\exp\left(-\sum_{m=1}^\infty
\frac{1}{m}\frac{\left(\frac{q_{2,1}}{q_{2,2}}\right)^m-
\left(\frac{p_{2,1}}{q_{2,2}}\right)^m}{1-\left(\frac{q_{3,2}}{q_{2,2}}\right)^m}w^m\right)
=
\exp\left(-\sum_{m=1}^\infty
\frac{1}{m}\frac{\left(\frac{q_{2,1}}{q_{2,2}}\right)^m-
\left(\frac{q_{2,1}}{p_{2,2}}\right)^m}{1-\left(\frac{q_{2,1}}{q_{1,1}}\right)^m}w^m\right).
\label{eqn:h21-1}
\end{eqnarray}
Here we used $|q_{j+1,j}/q_{j,j}|\neq 1$ $(1\leq j \leq 2)$ 
assumed in (\ref{assume:p,q1}).
From the compatibility of the two formulae for $h_{1,2}(w)$
in (\ref{eqn:h12-1})
[ or $h_{2,1}(w)$ in 
(\ref{eqn:h21-1})],
there are two possible restrictions 
for $q_{1,1}$, $q_{2,2}$, $q_{2,1}$,
and $q_{3,2}$,
\begin{eqnarray}
\mathrm{(i)}~~\frac{q_{2,1}}{q_{1,1}}=\frac{q_{3,2}}{q_{2,2}}~~~{\rm or}~~~
\mathrm{(ii)}~~\frac{q_{2,1}}{q_{1,1}}=\frac{q_{2,2}}{q_{3,2}}.
\label{choice:q}
\end{eqnarray}
 
First, we consider case 
$\mathrm{(i)}$~
$q_{2,1}/q_{1,1}=q_{3,2}/q_{2,2}$ in (\ref{choice:q}).
From the compatibility of the two formulae for $h_{1,2}(w)$
in (\ref{eqn:h12-1})
[ and $h_{2,1}(w)$
in (\ref{eqn:h21-1})],
we obtain
\begin{eqnarray}
&&
h_{1,2}(w)=\exp\left(-\sum_{m=1}^\infty
\frac{1}{m}\frac{[A_{1,2}(0)rm]_x}{[rm]_x}x^{-(A_{1,2}(0)+1)rm}
\left(\frac{q_{2,2}}{q_{1,1}}\right)^mw^m
\right),
\label{eqn:h12-3}
\\
&&
h_{2,1}(w)=\exp\left(-\sum_{m=1}^\infty
\frac{1}{m}\frac{[A_{1,2}(0)rm]_x}{[rm]_x}x^{(A_{1,2}(0)+1)rm}
\left(\frac{q_{1,1}}{q_{2,2}}\right)^mw^m
\right).
\label{eqn:h21-3}
\end{eqnarray}
We used (\ref{eqn:p,q2}) and $q_{2,1}/q_{1,1}=q_{3,2}/q_{2,2}=x^{2r}$.
From $h_{1,2}(w)=h_{2,1}(w)$ assumed in (\ref{assume:varphiSS}), we obtain
\begin{eqnarray}
\frac{q_{2,2}}{q_{1,1}}=x^{(A_{1,2}(0)+1)r}.
\label{eqn:p,q3}
\end{eqnarray}
Considering (\ref{eqn:p,q2}) and (\ref{eqn:p,q3}), we obtain
the first half of (\ref{eqn:p,q1}).
From (\ref{eqn:h12-3}), (\ref{eqn:h21-3}), and
(\ref{eqn:p,q3}), we obtain
\begin{eqnarray}
h_{1,2}(w)=h_{2,1}(w)=
\exp\left(-\sum_{m=1}^\infty
\frac{1}{m}\frac{[A_{1,2}(0)rm]_x}{[rm]_x}w^m
\right).
\label{sym:h}
\end{eqnarray}
Considering (\ref{normal3}) and (\ref{def:h}), we obtain the second half of (\ref{eqn:p,q1}).

Next, we consider case 
$\mathrm{(ii)}$~$
q_{2,1}/q_{1,1}=q_{2,2}/q_{3,2}$ in (\ref{choice:q}).
From the compatibility of the two formulae for $h_{1,2}(w)$
in (\ref{eqn:h12-1})
[and $h_{2,1}(w)$ in (\ref{eqn:h21-1})],
we obtain
\begin{eqnarray}
\left(\frac{p_{2,2}}{q_{1,1}}\right)^m+\left(\frac{q_{2,2}}{p_{2,1}}\right)^m=
\left(\frac{q_{2,2}}{q_{2,1}}\right)^m+\left(\frac{q_{2,2}}{q_{1,1}}\right)^m~~(m \neq 0).
\label{eqn:p,q4}
\end{eqnarray}
From (\ref{eqn:p,q4}) for $m=1,2$, we obtain 
$p_{2,2}/p_{2,1}=q_{2,2}/q_{2,1}$.
Combining (\ref{eqn:p,q4}) for $m=1$ 
and $p_{2,2}/p_{2,1}=q_{2,2}/q_{2,1}$, we obtain
$p_{2,1}=q_{2,1}$ or $p_{2,1}=q_{1,1}$.
For the case of $p_{2,1}=q_{2,1}$, we obtain $A_{1,2}(0)=0$ from (\ref{eqn:p,q2}).
For the case of $p_{2,1}=q_{1,1}$, we obtain $A_{1,2}(0)=1$ from (\ref{eqn:p,q2}).
$A_{1,2}(0)=0$ and $A_{1,2}(0)=1$ contradict with $-1<A_{1,2}(0)<0$ assumed in (\ref{assume:p,q1}).
Hence, the case $\mathrm{(ii)}~q_{2,1}/q_{1,1}=q_{2,2}/q_{3,2}$ is impossible.
~~~$\Box$

\begin{lem}
\label{lem:II-9}~
The relations (\ref{eqn:p,q5}) and (\ref{eqn:p,q5.5}) hold
if (\ref{assume:mutual-locality}), (\ref{assume:p,q1}), and (\ref{assume:start1}) are satisfied, where
\begin{eqnarray}
&&
\begin{array}{c}
{\displaystyle \frac{p_{j,j}}{q_{j,j}}=x^{A_{j,j}(0)r},~~~
\frac{p_{j+1,j}}{q_{j,j}}=x^{(2-A_{j,j}(0))r}},
\\
{\displaystyle s_j(m)s_j(-m)=-\frac{
\left[\frac{1}{2}A_{j,j}(0)rm\right]_x
\left[(2-A_{j,j}(0))rm\right]_x}
{\left[\frac{1}{2}(2-A_{j,j}(0))r m\right]_x [rm]_x}}
\end{array}
~~~(m>0,  A_{j,j}(0)\neq 1),
\label{eqn:p,q5}
\end{eqnarray}
\begin{eqnarray}
&&{\displaystyle p_{j,j}=p_{j+1,j},~~~
s_j(m)s_j(-m)=-1}~~~(m>0, A_{j,j}(0)=1).
\label{eqn:p,q5.5}
\end{eqnarray}
\end{lem}
{\it Proof.}~
From Lemma \ref{lem:II-7}, we obtain (\ref{eqn:hii-1}).
From (\ref{normal3}) and (\ref{def:h}), the constant term of the Taylor expansion for $h_{j,j}(w)$ is $1$.
Comparing the Taylor expansions for both sides of (\ref{eqn:hii-1}),
we obtain
\begin{eqnarray}
\left(\frac{q_{j+1,j}}{q_{j,j}}\right)^{A_{j,j}(0)-1}
\frac{p_{j+1,j}}{p_{j,j}}=1~~~(1\leq j \leq 2).
\label{eqn:p,q6}
\end{eqnarray}
Upon the specialization (\ref{eqn:p,q6}), the compatibility condition of the equations in
(\ref{eqn:hii-1}) is
\begin{eqnarray}
(p_{j,j}-p_{j+1,j})(p_{j,j}p_{j+1,j}-q_{j,j}q_{j+1,j})=0~~~(1\leq j \leq 2).
\label{eqn:p,q7}
\end{eqnarray}

First, we study the case of $A_{j,j}(0)=1$.
We obtain $p_{j,j}=p_{j+1,j}$ from (\ref{eqn:p,q6}).
Solving (\ref{eqn:hii-1}) upon $p_{j,j}=p_{j+1,j}$, we obtain 
$h_{j,j}(w)=1-w$.
Considering (\ref{normal3}) and (\ref{def:h}), we obtain $s_j(m)A_{j,j}(m)s_j(-m)=-1$ $(m>0)$.
We obtain (\ref{eqn:p,q5.5}).

Next, we study the case of $A_{j,j}(0)\neq 1$.
We obtain $p_{j,j}\neq p_{j+1,j}$ from (\ref{assume:p,q1}) and (\ref{eqn:p,q6}).
We obtain $p_{j,j}p_{j+1,j}=q_{j,j}q_{j+1,j}$ from $p_{j,j} \neq p_{j+1,j}$ and (\ref{eqn:p,q7}).
Combining $p_{j,j}p_{j+1,j}=q_{j,j}q_{j+1,j}$ and (\ref{eqn:p,q6}), we obtain the first part of (\ref{eqn:p,q5}).
Solving (\ref{eqn:hii-1}), we obtain
\begin{eqnarray}
h_{j,j}(w)=\exp\left(
-\sum_{m=1}^\infty \frac{1}{m}
{\displaystyle
\frac{\left[\frac{1}{2}A_{j,j}(0)rm\right]_x\left[(2-A_{j,j}(0))rm\right]_x}{
\left[\frac{1}{2}(2-A_{j,j}(0))r m\right]_x [rm]_x}}w^m
\right).\nonumber
\end{eqnarray}
We used 
$|q_{j+1,j}/q_{j,j}|\neq 1$ $(1\leq j \leq 2)$ in (\ref{assume:p,q1})
and
$q_{j+1,j}/q_{j,j}=x^{2r}$ $(1\leq j \leq 2)$ from (\ref{eqn:p,q1}).
Considering (\ref{normal3}) and (\ref{def:h}), we obtain the second part of (\ref{eqn:p,q5}).
~~$\Box$

\begin{prop}\label{prop:II-10}~
The relations (\ref{assume:mutual-locality}), (\ref{assume:q-difference}), and (\ref{assume:varphiSS}) hold
if the parameters $p_{i,j}, q_{i,j}$, $A_{i,j}(m)$, $s_i(m)$,
$g_i$, and $\lambda_{i,j}(m)$ are determined by (\ref{eqn:lambda(0)}), (\ref{eqn:lambda(m)}), (\ref{assume:start2}), 
(\ref{eqn:p,q1}), (\ref{eqn:p,q5}), (\ref{eqn:p,q5.5}), $p_{1,2}=q_{1,2}$, and $p_{3,1}=q_{3,1}$.
\end{prop}
{\it Proof.}~The proof is divided into Cases 1, 2, and 3 
classified in (\ref{def:4cases}) according to the values of $A_{1,1}(0)$ and $A_{2,2}(0)$.
For instance, for {Case1}, we will obtain the formulae
(\ref{CaseI-II-III:p,q,s(m)}), (\ref{CaseI-II-III:B}), (\ref{CaseI:A(0),g}), (\ref{CaseI:p,q}), (\ref{CaseI:A(m)-s(m)}), 
(\ref{CaseI:lambda(0)}), and (\ref{CaseI:lambda(m)}) in Theorem \ref{thm:II-1} from
the conditions (\ref{eqn:lambda(0)}), (\ref{eqn:lambda(m)}), (\ref{assume:start2}), (\ref{eqn:p,q1}), (\ref{eqn:p,q5}), (\ref{eqn:p,q5.5}), 
$p_{1,2}=q_{1,2}$, and $p_{3,1}=q_{3,1}$.
As a by-product of this calculation, we will show that
there is no indeterminacy in the free field realization
except for (\ref{def:shift}) and (\ref{def:scaling-boson}) upon the conditions (\ref{assume:p,q0}) and (\ref{assume:p,q1}).
Because $\lambda_{i,j}(m)$ are determined by (\ref{eqn:lambda(0)}) and (\ref{eqn:lambda(m)}),
the mutual locality (\ref{assume:mutual-locality}) holds from Lemma \ref{lem:II-5}.
The relations (\ref{assume:p,q0}) and (\ref{assume:start1}) are obtained by direct calculation using 
(\ref{CaseI:A(0),g}), (\ref{CaseI:p,q}), (\ref{CaseI:lambda(0)}), and (\ref{CaseI:lambda(m)}).
Hence, commutativity (\ref{assume:q-difference}) holds 
with (\ref{CaseI-II-III:B}) because of Lemma \ref{lem:II-6}.
Symmetry (\ref{assume:varphiSS}) holds because of (\ref{sym:h}).
Other cases are shown in the same way.
\\
$\bullet$~{Case 1} : From $A_{1,1}(0)\neq 1$ and (\ref{def:2parameters}), we set $A_{1,2}(0)=A_{2,1}(0)=(1-r)/r$.
From (\ref{eqn:p,q1}) and (\ref{eqn:p,q5}) we obtain
$p_{2,1}=sx^{2(1+A_{1,2}(0))r}=sx^{(2-A_{1,1}(0))r}$ and $p_{2,2}=sx^{(1-A_{1,2}(0))r}=sx^{(1+A_{1,2}(0)+A_{2,2}(0))r}$. 
Hence, we obtain $A_{1,1}(0)=A_{2,2}(0)=2(r-1)/r$ 
in (\ref{CaseI:A(0),g}).
Using the first half of
(\ref{CaseI:A(0),g}), (\ref{eqn:p,q1}), and (\ref{eqn:p,q5})
we obtain
$q_{i,i}=sx^{i-1}$, $q_{i+1,i}=sx^{2r+i-1}$, $p_{i,i}=sx^{2r+i-3}$, and $p_{i+1,i}=sx^{i+1}$ $(1\leq i \leq 2)$,
where $s=q_{1,1}$.
Upon the specialization $s=1$, we obtain (\ref{CaseI:p,q}).
Considering generic $s=q_{1,1}$, we obtain the second half of (\ref{def:shift}) for Case 1.
Using the first half of (\ref{CaseI:A(0),g}), (\ref{CaseI:p,q}), and (\ref{assume:start2}), we obtain $g_1=g_2=g_3$ in the second half of (\ref{CaseI:A(0),g}).
From (\ref{CaseI:A(0),g}), (\ref{eqn:p,q1}), and (\ref{eqn:p,q5}) we obtain
\begin{eqnarray}
s_i(m)A_{i,j}(m)s_j(-m)=\left\{\begin{array}{cc}
{\displaystyle -\frac{[(r-1)m]_x [2m]_x}{[rm]_x[m]_x}}&(1\leq i=j \leq 2),
\\
{\displaystyle \frac{[(r-1)m]_x}{[rm]_x}}&(1\leq i \neq j \leq 2)
\end{array}
\right.~~~(m>0).
\nonumber
\end{eqnarray}
Hence, we obtain
\begin{eqnarray}
s_j(-m)=-\frac{[(r-1)m]_x [2m]_x}{[rm]_x[m]_x}\frac{1}{s_j(m)}~~~(m>0, 1\leq j \leq 2),
\nonumber
\\
A_{i,j}(m)=A_{j,i}(-m)=-\frac{[m]_x}{[2m]_x}\frac{s_j(m)}{s_i(m)}~~~(m>0, 1\leq i \neq j \leq 2).
\nonumber
\end{eqnarray}
Thus, setting $s_j(m)=1$~$(m>0, 1\leq j \leq 2)$ provides 
(\ref{CaseI:A(m)-s(m)}).
Setting $s_j(m)=\alpha_j(m)\neq 0$ $(m>0, 1\leq j \leq 2)$ provides the scaling of the free field (\ref{def:scaling-boson})
for Case1.
Solving the linear equations (\ref{eqn:lambda(0)}) and (\ref{eqn:lambda(m)}) for $\lambda_{i,j}(m)$
upon (\ref{CaseI-II-III:p,q,s(m)}), (\ref{CaseI:A(0),g}), (\ref{CaseI:p,q}), and (\ref{CaseI:A(m)-s(m)}),
we obtain (\ref{CaseI:lambda(0)}) and (\ref{CaseI:lambda(m)}).
The first half of (\ref{def:shift}) for Case 1 is obtained in the same way.
We obtained (\ref{CaseI-II-III:p,q,s(m)}),
(\ref{CaseI-II-III:B}) for Case 1, and (\ref{CaseI:A(0),g})--(\ref{CaseI:lambda(m)}) in Theorem \ref{thm:II-1}.
We obtained (\ref{def:shift}) and (\ref{def:scaling-boson}) for Case 1 from necessary conditions.
\\
$\bullet$ {Case 2} :
From $A_{1,1}(0)\neq 1$ and (\ref{def:2parameters}), we set $A_{1,2}(0)=A_{2,1}(0)=(1-r)/r$.
From (\ref{eqn:p,q1}) and (\ref{eqn:p,q5}), we obtain
$p_{2,1}=sx^{2(1+A_{1,2}(0))r}=sx^{(2-A_{1,1}(0))r}$. 
Hence, we obtain $A_{1,1}(0)=2(r-1)/r$ in (\ref{CaseII:A(0),g}).
Using the first half of (\ref{CaseII:A(0),g}), (\ref{eqn:p,q1}), and (\ref{eqn:p,q5}), we obtain 
$q_{i,i}=sx^{i-1}$, $q_{i+1,i}=sx^{2r+i-1}$, $p_{i,i}=sx^{2r+i-3}$, and $p_{i+1,i}=sx^{-2r+5+(2r-3)i}$ $(1\leq i \leq 2)$ where 
$s=q_{1,1}$.
Upon the specialization $s=1$, we obtain (\ref{CaseII:p,q}).
Considering generic $s=q_{1,1}$, we obtain the second half of (\ref{def:shift}) for Case 2.
Using the first half of (\ref{CaseII:A(0),g}), (\ref{CaseII:p,q}), and (\ref{assume:start2}), we obtain $g_1=g_2$ and 
$g_3=[r-1]_x g_1$ in the second half of (\ref{CaseII:A(0),g}).
From (\ref{CaseII:A(0),g}), (\ref{eqn:p,q1}), (\ref{eqn:p,q5}), and (\ref{eqn:p,q5.5}), we obtain
\begin{eqnarray}
s_i(m)A_{i,j}(m)s_j(-m)=\left\{\begin{array}{cc}
{\displaystyle -\frac{[(r-1)m]_x [2m]_x}{[rm]_x[m]_x}}&(i=j=1),
\\
-1&(i=j=2),
\\
{\displaystyle \frac{[(r-1)m]_x}{[rm]_x}}&(1\leq i \neq j \leq 2)
\end{array}
\right.~~~(m>0).
\nonumber
\end{eqnarray}
Hence, we obtain
\begin{eqnarray}
&&
s_1(-m)=-\frac{[(r-1)m]_x [2m]_x}{[rm]_x[m]_x}\frac{1}{s_1(m)},
~~~s_2(-m)=\frac{1}{s_2(m)}~~~(m>0),
\nonumber
\\
&&
A_{1,2}(m)=-\frac{[(r-1)m]_x}{[rm]_x}\frac{s_2(m)}{s_1(m)},
~~~
A_{2,1}(m)=-\frac{[m]_x}{[2m]_x}\frac{s_1(m)}{s_2(m)}~~~
(m>0).
\nonumber
\end{eqnarray}
Thus, setting $s_j(m)=1$~$(m>0, 1\leq j \leq 2)$ provides 
(\ref{CaseII:A(m)-s(m)}).
Setting $s_j(m)=\alpha_j(m)\neq 0$ $(m>0, 1\leq j \leq 2)$ provides the scaling of the free field (\ref{def:scaling-boson}) for Case 2.
Solving the linear equations (\ref{eqn:lambda(0)}) and (\ref{eqn:lambda(m)}) for $\lambda_{i,j}(m)$ upon 
(\ref{CaseI-II-III:p,q,s(m)}), (\ref{CaseII:A(0),g}), (\ref{CaseII:p,q}), and (\ref{CaseII:A(m)-s(m)}),
we obtain (\ref{CaseII:lambda(0)}) and (\ref{CaseII:lambda(m)}).
The second half of (\ref{def:shift}) for Case 2 is obtained in the same way.
We obtained 
(\ref{CaseI-II-III:p,q,s(m)}), (\ref{CaseI-II-III:B}) for Case 2
and (\ref{CaseII:A(0),g})--(\ref{CaseII:lambda(m)}) in Theorem \ref{thm:II-1}.
We obtained (\ref{def:shift}) and (\ref{def:scaling-boson}) for Case 2 from necessary conditions.
\\
$\bullet$ {Case 3} :
From $A_{1,1}(0)=1$ and (\ref{def:2parameters}), we set $A_{1,2}(0)=A_{2,1}(0)=-1/r$.
Using the first half of (\ref{CaseIII:A(0),g}), (\ref{eqn:p,q1}), and (\ref{eqn:p,q5.5}), we obtain 
$q_{i,i}=sx^{(r-1)(i-1)}$, $q_{i+1,i}=sx^{r+1+(r-1)i}$, $p_{i,i}=p_{i+1,i}=sx^{3r-5+(-r+3)i}$ $(1\leq i \leq 2)$,
where $s=q_{1,1}$.
Upon the specialization $s=1$, we obtain (\ref{CaseIII:p,q}).
Considering generic $s=q_{1,1}$, we obtain the second half of (\ref{def:shift}) for Case 3.
Using the first half of (\ref{CaseIII:A(0),g}), (\ref{CaseIII:p,q}), and (\ref{assume:start2}), we obtain $g_2=[r-1]_x g_1$ and 
$g_3=g_1$ in the second half of (\ref{CaseIII:A(0),g}).
From (\ref{CaseIII:A(0),g}), (\ref{eqn:p,q1}), and (\ref{eqn:p,q5.5}), 
we obtain
\begin{eqnarray}
s_i(m)A_{i,j}(m)s_j(-m)=\left\{\begin{array}{cc}
-1 &(1\leq i=j \leq 2),
\\
{\displaystyle \frac{[m]_x}{[rm]_x}}&(1\leq i \neq j \leq 2)
\end{array}
\right.~~~(m>0).
\nonumber%\label{eqn:sAs5}
\end{eqnarray}
Hence, we obtain
\begin{eqnarray}
s_j(-m)=-\frac{1}{s_j(m)}~~(m>0, 1\leq j \leq 2),
~~A_{i,j}(m)=A_{j,i}(-m)=-\frac{[m]_x}{[rm]_x}\frac{s_j(m)}{s_i(m)}~~(m>0, 1\leq i \neq j \leq 2).
\nonumber
\end{eqnarray}
Thus, setting $s_j(m)=1$~$(m>0, 1\leq j \leq 2)$ provides 
(\ref{CaseIII:A(m)-s(m)}).
Setting $s_j(m)=\alpha_j(m)\neq 0$ $(m>0, 1\leq j \leq 2)$ provides the scaling of the free field (\ref{def:scaling-boson}) for Case 3.
Solving the linear equations (\ref{eqn:lambda(0)}) and (\ref{eqn:lambda(m)}) for $\lambda_{i,j}(m)$
upon (\ref{CaseI-II-III:p,q,s(m)}), (\ref{CaseIII:A(0),g}), (\ref{CaseIII:p,q}), 
and (\ref{CaseIII:A(m)-s(m)}),
we obtain (\ref{CaseIII:lambda(0)}) and (\ref{CaseIII:lambda(m)}).
The first half of (\ref{def:shift}) is derived in the same way.
We obtained (\ref{CaseI-II-III:p,q,s(m)}), (\ref{CaseI-II-III:B}) for Case 3 and 
(\ref{CaseIII:A(0),g})--(\ref{CaseIII:lambda(m)}) in Theorem \ref{thm:II-1}.
We obtained (\ref{def:shift}) and (\ref{def:scaling-boson}) for
Case 3 from necessary conditions.
~~$\Box$

In the proof of Proposition \ref{prop:II-10}
we proved Theorem \ref{thm:II-1} at the same time.
As a by-product, we proved 
that there in no determinacy in the realization
except for (\ref{def:arrange-suffix}), (\ref{def:shift}), and (\ref{def:scaling-boson}) upon the conditions 
(\ref{assume:p,q0}) and  (\ref{assume:p,q1}).

~\\
{\it Proof of Proposition \ref{prop:II-3}.}~
Using $h_{k,l}(w)$ in (\ref{def:h}) we obtain
\begin{eqnarray}
S_k(w_1)S_l(w_2)=\left(\frac{w_1}{w_2}\right)^{A_{k,l}(0)}\frac{h_{k,l}\left(\frac{w_2}{w_1}\right)}{
h_{l,k}\left(\frac{w_1}{w_2}\right)}
S_l(w_2)S_k(w_1)~~~(1\leq k, l \leq 2).
\nonumber
\end{eqnarray}
Using (\ref{eqn:p,q1}), (\ref{eqn:p,q5}) and (\ref{eqn:p,q5.5}) 
we obtain (\ref{CaseI:S1S2}), (\ref{CaseI:SjSj}), (\ref{CaseII:SjSj}), (\ref{CaseIII:S1S2}), and (\ref{CaseIII:SjSj}).
~~$\Box$\\

\section{Quadratic relations}
\label{Section3}
In this section we introduce the higher $W$-currents $T_i(z)$
and obtain a set of quadratic relations of $T_i(z)$ 
for the deformed
$W$-superalgebra ${\cal W}_{q,t}(\mathfrak{sl}(2|1))$.

\subsection{Quadratic relations}

We define the functions $\Delta_i(z)$ $(i=0,1,2, \cdots)$ as
\begin{eqnarray}
\Delta_i(z)=\frac{(1-x^{2r-i}z)(1-x^{-2r+i}z)}{(1-x^i z)(1-x^{-i}z)}.
\label{def:Delta}
\end{eqnarray}
We have
\begin{eqnarray}
&&
\Delta_i(z)-\Delta_i({z}^{-1})=\frac{[r]_x[r-i]_x}{[i]_x}(x-x^{-1})(\delta(x^{-i}z)-\delta(x^iz))~~(i=1,2,3,\cdots).
\nonumber
%\label{eqn:delta}
\end{eqnarray}

We define the structure functions $f_{i,j}(z)$ $(i, j=0,1,2,\cdots)$ as
\begin{eqnarray}
f_{i,j}(z)=\exp\left(-\sum_{m=1}^\infty \frac{1}{m}
\frac{[(r-1)m]_x[rm]_x
[{\rm Min}(i,j)m]_x[(s-{\rm Max}(i,j))m]_x
}{[m]_x[sm]_x}(x-x^{-1})^2
z^m
\right),
\label{CaseI-II-III:fij}
\end{eqnarray}
where
\begin{eqnarray}
s=\left\{
\begin{array}{cc}
3&for~Case~1,\\
r+1& for~Cases~2~and~3.
\end{array}
\right.\nonumber
\end{eqnarray}
The ratio of the structure function
\begin{eqnarray}
\frac{f_{1,1}(z^{-1})}{f_{1,1}(z)}=-z
\frac{\Theta_{x^{2s}}(x^2z, x^{2s-2r}z, x^{2s+2r-2}z)}{
\Theta_{x^{2s}}(x^2z^{-1},x^{2s-2r}z^{-1}, x^{2s+2r-2}z^{-1})}
\nonumber
\end{eqnarray}
coincides with those of (\ref{CaseI-II-III:vertex}).

We introduce the higher $W$-currents $T_i(z)$ and give the quadratic relations.
From now on, we set $g_1=1$
for all cases, but this is not an essential limitation.
Hereafter, we use the abbreviations
\begin{eqnarray}
c(r, x)=[r]_x[r-1]_x(x-x^{-1}),~~~d_0(r,x)=1,~~d_j(r, x)=\prod_{l=1}^j \frac{[r-l]_x}{[l]_x}~~~(j \geq 1).
\label{def:abbreviation}
\end{eqnarray}

For Case 1, we set the $W$-currents as
\begin{align}
&
T_1(z)=\Lambda_1(z)+\Lambda_2(z)+\Lambda_3(z),\notag
\\
&
T_2(z)=:\Lambda_1(x^{-1}z)\Lambda_2(xz):
+:\Lambda_1(x^{-1}z)\Lambda_3(xz):+:\Lambda_2(x^{-1}z)\Lambda_3(xz):.
\notag
%\label{CaseI:Ti(z)}
\end{align}
They satisfy the quadratic relations of ${\cal W}_{q,t}(\mathfrak{sl}(3))$\cite{Awata-Kubo-Odake-Shiraishi, Feigin-Frenkel}
\begin{eqnarray}
&&
f_{1,1}\left(\frac{z_2}{z_1}\right)T_1(z_1)T_1(z_2)-f_{1,1}\left(\frac{z_1}{z_2}\right)T_1(z_2)T_1(z_1)
=c(r,x)\left(
\delta\left(\frac{x^{-2}z_2}{z_1}\right)T_2(x^{-1}z_2)
-\delta\left(\frac{x^2z_2}{z_1}\right)T_2(x z_2)
\right),\nonumber\\
&&
f_{1,2}\left(\frac{z_2}{z_1}\right)T_1(z_1)T_2(z_2)-f_{2,1}\left(\frac{z_1}{z_2}\right)T_2(z_2)T_1(z_1)
=c(r,x)\left(
\delta\left(\frac{x^{-3}z_2}{z_1}\right)
-\delta\left(\frac{x^3z_2}{z_1}\right)
\right),\label{CaseI:quadratic}
\\
&&
f_{2,2}\left(\frac{z_2}{z_1}\right)T_2(z_1)T_2(z_2)-f_{2,2}\left(\frac{z_1}{z_2}\right)T_2(z_2)T_2(z_1)
=c(r,x)\left(
\delta\left(\frac{x^{-2}z_2}{z_1}\right)T_1(x^{-1}z_2)
-\delta\left(\frac{x^2z_2}{z_1}\right)T_1(x z_2)
\right),\nonumber
\end{eqnarray}
with $f_{1,1}(z)=f_{2,2}(z)$ and $f_{1,2}(z)=f_{2,1}(z)$.
Here, we omit the proof. 
The deformed $W$-algebra ${\cal W}_{q, t}(\mathfrak{sl}(3))$
is the associative algebra over ${\mathbf C}$ with the generators $T_i[m]$ $(m \in {\mathbf Z}, i=1,2)$
and defining relations (\ref{CaseI:quadratic}).
Here we set $T_i(z)=\sum_{m \in {\mathbf Z}}T_i[m]z^{-m}$ $(i=1,2)$.
The parameters $q$ and $t$ in ${\cal W}_{q,t}(\mathfrak{sl}(3))$
are given as $q=x^{2r}$ and $t=x^{2(r-1)}$.

For Case 2, we introduce the $W$-currents $T_i(z)$ $(i=1,2,3,\cdots)$ as
\begin{align}
T_1(z)&=\Lambda_1(z)+\Lambda_2(z)+
d_1(r,x)\Lambda_3(z),
\notag
\\
T_2(z)&=:\Lambda_1(x^{-1}z)\Lambda_2(xz):+d_1(r,x):\Lambda_1(x^{-1}z)\Lambda_3(xz):+
d_1(r,x):\Lambda_2(x^{-1}z)\Lambda_3(xz):
\notag
\\
&+d_2(r,x):\Lambda_3(x^{-1}z)\Lambda_3(xz):,\notag
\\
T_i(z)
&=d_{i-2}(r,x)
:\Lambda_1(x^{-i+1}z) \Lambda_2(x^{-i+3}z)
\prod_{j=1}^{i-2} \Lambda_3(x^{-i+2j+3}z):
\notag\\
&+d_{i-1}(r,x)
:\Lambda_1(x^{-i+1}z) \prod_{j=1}^{i-1} \Lambda_3(x^{-i+2j+1}z):
+d_{i-1}(r,x)
:\Lambda_2(x^{-i+1}z)
\prod_{j=1}^{i-1} \Lambda_3(x^{-i+2j+1}z):
\notag\\
&+d_i(r,x):\prod_{j=1}^{i} \Lambda_3(x^{-i+2j-1}z):~~(i=3,4,5,\cdots).
\label{CaseII:Ti(z)}
\end{align}
We use $d_j(r, x)$ defined in (\ref{def:abbreviation}).
We have $T_i(z)\neq 1$ $(i=1,2,3,\cdots)$ and $T_i(z) \neq T_j(z)$ $(i \neq j)$.
The definition of $T_i(z)$ for Case 3 is similar as those for Case 2.
See the definition of $T_i(z)$ for Case 3 in (\ref{CaseIII:Ti(z)}).
On the other hand, the definitions of $T_i(z)$ for Cases 2 and 3 are different from those of Case 1.
See discussion in Section \ref{Section4}.

The following is {\bf the main theorem} of this paper which holds for Cases 2 and 3.

\begin{thm}
\label{thm:III-1}~
For the deformed $W$-superalgebra ${\cal W}_{q,t}(\mathfrak{sl}(2|1))$
the $W$-currents $T_i(z)$
satisfy the set of quadratic relations
\begin{align}
&
f_{i,j}\left(\frac{z_2}{z_1}\right)T_i(z_1)T_j(z_2)-f_{j,i}\left(\frac{z_1}{z_2}\right)T_j(z_2)T_i(z_1)\notag\\
=&~c(r,x)\sum_{k=1}^i \prod_{l=1}^{k-1}\Delta_1(x^{2l+1})
\left(
\delta\left(\frac{x^{-j+i-2k}z_2}{z_1}\right)f_{i-k, j+k}(x^{j-i})T_{i-k}(x^{k}z_1)T_{j+k}(x^{-k}z_2)
\right.\notag
\\
-&\left.\delta\left(\frac{x^{j-i+2k}z_2}{z_1}\right)f_{i-k, j+k}(x^{-j+i})T_{i-k}(x^{-k}z_1)T_{j+k}(x^kz_2)
\right)~~~(j \geq i \geq 1).
\label{CaseII-III:quadratic}
\end{align}
Here we use $f_{i,j}(z)$ in (\ref{CaseI-II-III:fij}).
$T_0(z)$ in the right hand side is understood as $T_0(z)=1$.
\end{thm}

In view of Theorem \ref{thm:III-1} we arrive at the following definition. 
\begin{dfn}
Set $T_i(z)=\sum_{m \in {\mathbf Z}}T_i[m]z^{-m}$ $(i=1,2,3,\cdots)$.
The deformed $W$-superalgebra ${\cal W}_{q, t}(\mathfrak{sl}(2|1))$
is the associative algebra over ${\mathbf C}$ with the generators $T_i[m]$ $(m \in {\mathbf Z}, i=1,2,3,\cdots)$
and defining relations (\ref{CaseII-III:quadratic}).
\end{dfn}

\subsection{Proof of Theorem \ref{thm:III-1}}

\begin{lem}
\label{lem:III-2}
For Cases 2 and 3, $\Delta_i(z)$ and $f_{i,j}(z)$ satisfy the following fusion relations.
\begin{eqnarray}
&&
f_{i,j}(z)=f_{j,i}(z)=\prod_{k=1}^i f_{1,j}(z^{-i-1+2k}z)~~~(1\leq i \leq j),
\label{eqn:fusion0}
\\
&&
f_{1,i}(z)=
\left(
\prod_{k=1}^{i-1}\Delta_1(x^{-i+2k}z)
\right)^{-1}
\prod_{k=1}^i f_{1,1}(x^{-i-1+2k}z)~~~
(i \geq 2),
\label{eqn:fusion5}
\\
&&
f_{1,i}(z)f_{j,i}(x^{\pm (j+1)}z)=\left\{\begin{array}{cc}
f_{j+1,i}(x^{\pm j}z)\Delta_1(x^{\pm i}w)&(1\leq i \leq j),\\
f_{j+1,i}(x^{\pm j}z)&(1\leq j<i),
\end{array}\right.
\label{eqn:fusion1}\\
&&
f_{1,i}(z)f_{1,j}(x^{\pm (i+j)}z)=
f_{1,i+j}(x^{\pm j}z)\Delta_1(x^{\pm i}z)~~(i,j \geq 1),
\label{eqn:fusion2}\\
&&
f_{1,i}(z)f_{1,j}(x^{\pm (i-j-2k)}z)=f_{1,i-k}(x^{\mp k}z)f_{1,j+k}(x^{\pm (i-j-k)}z)~~(i,j,i-k,j+k \geq 1),
\label{eqn:fusion3}
\\
&&
\Delta_{i+1}(z)=
\left(\prod_{k=1}^{i-1} \Delta_{1}(x^{-i+2k}z)
\right)^{-1}
\prod_{k=1}^{i}\Delta_2(x^{-i-1+2k}z)~~~
(i \geq 2).
\label{eqn:fusion4}
\end{eqnarray}
\end{lem}
{\it Proof.}~
We obtain (\ref{eqn:fusion0}) and (\ref{eqn:fusion4})
by straightforward calculation from the definitions.
We show (\ref{eqn:fusion5}) here.
From definitions (\ref{def:Delta}) and (\ref{CaseI-II-III:fij}), we have
\begin{align}
&
\left(\prod_{k=1}^{i-1}\Delta_1(x^{-i+2k}z)\right)^{-1}
\prod_{k=1}^i f_{1,1}(x^{-i-1+2k}z)\nonumber
\\
=&
\exp\left(
-\sum_{m=1}^\infty
\frac{1}{m}
\frac{[rm]_x[(r-1)m]_x}{[(r+1)m]_x}
(x-x^{-1})^2
\left([rm]_x\sum_{k=1}^i x^{(-i+2k-1)m}-
[(r+1)m]_x\sum_{k=1}^{i-1}
x^{(-i+2k)m}
\right)z^m
\right).\nonumber
\end{align}
Using the relation
\begin{eqnarray}
[rm]_x\sum_{k=1}^i x^{(-i+2k-1)m}-
[(r+1)m]_x\sum_{k=1}^{i-1}
x^{(-i+2k)m}=[(r+1-i)m]_x,
\nonumber
\end{eqnarray}
we have $f_{1,i}(z)$.
Using (\ref{eqn:fusion0}) and (\ref{eqn:fusion5}),
we obtain the relations (\ref{eqn:fusion1}), (\ref{eqn:fusion2}), and (\ref{eqn:fusion3}).~~~
$\Box$

\begin{prop}
\label{prop:III-3}
$\Lambda_i(z)$'s satisfy
\begin{eqnarray}
&&
~f_{1,1}\left(\frac{z_2}{z_1}\right)
\Lambda_k(z_1)\Lambda_l(z_2)=
\Delta_1\left(\frac{x^{-1}z_2}{z_1}\right)
:\Lambda_k(z_1)\Lambda_l(z_2):,
\nonumber
\\
&&
f_{1,1}\left(\frac{z_2}{z_1}\right)
\Lambda_l(z_1)\Lambda_k(z_2)=
\Delta_1\left(\frac{x z_2}{z_1}\right)
:\Lambda_l(z_1)\Lambda_k(z_2):
~~~(1\leq k<l\leq 3),
\label{CaseI-II-III:lambda_{kl}}
\\
&&
f_{1,1}\left(\frac{z_2}{z_1}\right)
\Lambda_1(z_1)\Lambda_1(z_2)=
:\Lambda_1(z_1)\Lambda_1(z_2):,
\nonumber
\\
&&
f_{1,1}\left(\frac{z_2}{z_1}\right)
\Lambda_i(z_1)\Lambda_i(z_2)=
\Delta_2\left(\frac{s_i z_2}{z_1}\right):
\Lambda_i(z_1)\Lambda_i(z_2):
(2 \leq i \leq 3),
\label{CaseI-II-III:lambda_{ii}}
\end{eqnarray}
where
\begin{eqnarray}
s_2=\left\{\begin{array}{cc}
0&for~Cases~1~and~2,\\
1&for~Case~3,
\end{array}\right.~~~
s_3=\left\{
\begin{array}{cc}
0&for~Cases~1~and~3,\\
1&for~Case~2.
\end{array}
\right.
\nonumber
\end{eqnarray}
Here we use $f_{1,1}(z)$ in (\ref{CaseI-II-III:fij}).
\end{prop}
{\it Proof.}~
Substituting 
(\ref{CaseI:lambda(m)}), (\ref{CaseII:lambda(m)}), and (\ref{CaseIII:lambda(m)}) into 
$\varphi_{\Lambda_k, \Lambda_l}(z_1,z_2)$ in (\ref{normal4}) we obtain (\ref{CaseI-II-III:lambda_{kl}}) and
(\ref{CaseI-II-III:lambda_{ii}}).~~~
$\Box$

~\\
{\it Proof of Proposition \ref{prop:II-2}.}
\label{proof:prop:II-2}~
Using $\varphi_{\Lambda_k, \Lambda_l}(z_1,z_2)$ in (\ref{normal4}) we obtain
\begin{eqnarray}
\Lambda_k(z_1)\Lambda_l(z_2)=\frac{\varphi_{\Lambda_k, \Lambda_l}\left(z_1, z_2\right)}{
\varphi_{\Lambda_l, \Lambda_k}\left(z_2, z_1\right)}
\Lambda_l(z_2)\Lambda_k(z_1)~~~(1\leq k, l \leq 3).
\nonumber
\end{eqnarray}
Using the explicit formulae of $\varphi_{\Lambda_k, \Lambda_l}(z_1,z_2)$,
we obtain (\ref{CaseI-II-III:vertex}).
~~$\Box$\\

Here, we set the higher $W$-currents $T_i(z)$ for Case 3.
To avoid confusion, we write $\Lambda_i(z)$ for Case 2 
(resp. Case 3) as 
$\Lambda_i^{\rm II}(z)$ (resp. $\Lambda_i^{\rm III}(z)$).
We write $T_i(z)$ for Case 2 (resp. Case 3) as $T_i^{\rm II}(z)$ (resp. $T_i^{\rm III}(z)$).
From Proposition \ref{prop:III-3},
$\Lambda_i^{\rm II}(z)$ and $\Lambda_i^{\rm III}(z)$ 
satisfy the same relations (\ref{CaseI-II-III:lambda_{kl}}), (\ref{CaseI-II-III:lambda_{ii}}), where $i \to i'$ is the permutation $(1,2,3) \to (3,1,2).$
If the expression of the $W$-currents $T_i^{\rm II}(z)$ in (\ref{CaseII:Ti(z)}) is abbreviated as
\begin{eqnarray}
T_i^{\rm II}(z)=P_i(\Lambda_1^{\rm II}(z), \Lambda_2^{\rm II}(z), \Lambda_3^{\rm II}(z))~~~(i=1,2,3,\cdots),
\nonumber
\end{eqnarray}
the $W$-currents $T_i^{\rm III}(z)$ is defined as
\begin{eqnarray}
T_i^{\rm III}(z)=P_i(\Lambda_3^{\rm III}(z), \Lambda_1^{\rm III}(z), \Lambda_2^{\rm III}(z))
~~~(i=1,2,3,\cdots).\label{CaseIII:Ti(z)}
\end{eqnarray}
The algebras generated by $T_i(z)$ $(i=1,2,3,\cdots)$ are the same as for Case 2 and Case 3.

The following lemmas \ref{lem:III-4}, \ref{lem:III-5}, and \ref{lem:III-6} give special cases of (\ref{CaseII-III:quadratic}).
\begin{lem}
\label{lem:III-4}~
For Cases 2 and 3 we have the fusion relation of $T_i(z)$ as
\begin{eqnarray}
&&
\lim_{z_1 \to x^{\pm (i+j)}z_2}
\left(1-x^{\pm(i+j)}\frac{z_2}{z_1}\right)f_{i,j}\left(\frac{z_2}{z_1}\right)T_i(z_1)T_j(z_2)\nonumber
\\
&&=\mp c(r,x) \prod_{k=1}^{{\rm Min}(i,j)-1}\Delta_1(x^{2k+1})T_{i+j}(x^{\pm i}z_2)~~~(i,j \geq 1).
\label{eqn:fusion6}
\end{eqnarray}
\end{lem}
{\it Proof.}~
Summing up the relations 
(\ref{appendixA:1})--(\ref{appendixA:12}) in Appendix \ref{Appendix:fusion}
gives (\ref{eqn:fusion6}).~~$\Box$

\begin{lem}
\label{lem:III-5}~
For Cases 2 and 3 we have the exchange relation as 
meromorphic functions
\begin{eqnarray}
f_{i,j}\left(\frac{z_2}{z_1}\right)T_i(z_1)T_j(z_2)
=f_{j,i}\left(\frac{z_1}{z_2}\right)T_j(z_2)T_i(z_1)~~~(j \geq i \geq 1).
\label{eqn:TiTj-commute}
\end{eqnarray}
\end{lem}
{\it Proof.}~Using the commutation relations 
(\ref{CaseI-II-III:lambda_{kl}}) 
and (\ref{CaseI-II-III:lambda_{ii}}) repeatedly, 
(\ref{eqn:TiTj-commute}) is obtained except for poles in both sides.
~~$\Box$

\begin{lem}
\label{lem:III-6}~
For Cases 2 and 3 $T_i(z)$ satisfy the quadratic relation.
\begin{align}
&
f_{1,i}\left(\frac{z_2}{z_1}\right)T_1(z_1)T_i(z_2)-f_{i,1}\left(\frac{z_1}{z_2}\right)T_i(z_2)T_1(z_1)\nonumber\\
=&~c(r,x)\left(\delta\left(\frac{x^{-i-1}z_2}{z_1}\right)T_{i+1}(x^{-1}z_2)-\delta\left(\frac{x^{i+1}z_2}{z_1}\right)T_{i+1}(xz_2)\right)
~~(i \geq 1).
\label{def:quadratic'}
\end{align}
\end{lem}
{\it Proof}.~
Summing up the relations 
(\ref{appendixB:1})--(\ref{appendixB:7}) in Appendix \ref{Appendix:exchange} 
gives the quadratic relation (\ref{def:quadratic'}).~~
$\Box$

~\\
~{\it Proof of Theorem \ref{thm:III-1}}.~
We prove Theorem \ref{thm:III-1} by induction.
Lemma \ref{lem:III-6} is the basis of induction for the proof.
\color{black}
We define ${\rm LHS}_{i,j}$ and ${\rm RHS}_{i,j}(k)$ with $(1\leq k \leq i \leq j)$ as
\begin{align}
{\rm LHS}_{i,j}&=
f_{i,j}\left(\frac{z_2}{z_1}\right)T_i(z_1)T_j(z_2)-f_{j,i}\left(\frac{z_1}{z_2}\right)T_j(z_2)T_i(z_1),
\notag
\\
{\rm RHS}_{i,j}(k)&=
c(r,x)
\prod_{l=1}^{k-1}\Delta_1(x^{2l+1})
\left(\delta\left(\frac{x^{-j+i-2k}z_2}{z_1}\right)
f_{i-k,j+k}(x^{j-i})T_{i-k}(x^kz_1)T_{j+k}(x^{-k}z_2)\right.
\notag\\
&-\left.\delta
\left(\frac{x^{j-i+2k}z_2}{z_1}\right)f_{i-k,j+k}(x^{-j+i})T_{i-k}(x^{-k}z_1)T_{j+k}(x^{k}z_2)
\right)~~
(1\leq k \leq i-1),
\notag
\\
{\rm RHS}_{i,j}(i)&=
c(r,x)
\prod_{l=1}^{i-1}\Delta_1(x^{2l+1})
\left(\delta\left(\frac{x^{-j-i}z_2}{z_1}\right)T_{j+i}(x^{-i}z_2)-
\delta\left(\frac{x^{j+i}z_2}{z_1}\right)
T_{j+i}(x^{i}z_2)\right).
\notag
\end{align}
We prove the following relation by induction on 
$i$ $(1\leq i \leq j)$.
\begin{eqnarray}
{\rm LHS}_{i,j}=\sum_{k=1}^i {\rm RHS}_{i,j}(k).
\label{induction4}
\end{eqnarray}
The starting point of $i=1 \leq j$ was previously proven in 
Lemma \ref{lem:III-6}.
We assume that relation (\ref{induction4}) holds for some $i$ $(1\leq i<j)$, and we show ${\rm LHS}_{i+1,j}=\sum_{k=1}^{i+1} {\rm RHS}_{i+1,j}(k)$ from this assumption.
Multiplying ${\rm LHS}_{i,j}$ by 
$f_{1,i}\left(z_1/z_3\right)f_{1,j}\left(z_2/z_3\right)T_1(z_3)$ on the left and using the quadratic relation (\ref{induction4}) with $i=1$
along with fusion relation (\ref{eqn:fusion1}) gives 
\begin{align}
&f_{1,j}\left(\frac{z_2}{z_3}\right)f_{i,j}\left(\frac{z_2}{z_1}\right)f_{1,i}\left(\frac{z_1}{z_3}\right)T_1(z_3)T_i(z_1)T_j(z_2)
-f_{j,1}\left(\frac{z_3}{z_2}\right)f_{j,i}\left(\frac{z_1}{z_2}\right)T_j(z_2)
f_{1,i}\left(\frac{z_1}{z_3}\right)T_1(z_3)T_i(z_1)\notag\\
-&~c(r,x)\delta\left(\frac{x^{-j-1}z_2}{z_3}\right)
\Delta_1\left(\frac{x^{-i}z_1}{z_3}\right)f_{j+1,i}\left(\frac{x^{-j}z_1}{z_3}\right)T_{j+1}(x^jz_3)T_i(z_1)\notag
\\
+&~c(r,x)\delta\left(\frac{x^{j+1}z_2}{z_3}\right)
\Delta_1\left(\frac{x^{i}z_1}{z_3}\right)f_{j+1,i}\left(\frac{x^{j}z_1}{z_3}\right)T_{j+1}(x^{-j}z_3)T_i(z_1).
\label{proof1}
\end{align}
Taking the limit $z_3\to x^{-i-1}z_1$ of (\ref{proof1}) multiplied by $c(r,x)^{-1} \left(1-x^{-i-1}z_1/z_3\right)$ and using fusion relation (\ref{eqn:fusion6}) along with the relation $\lim_{z_3 \to x^{-i-1}z_1}\left(1-x^{-i-1}z_1/z_3\right)
\Delta_1\left(x^{-i}z_1/z_3\right)=c(r,x)$ gives
\begin{align}
&f_{1,j}\left(\frac{x^{i+1}z_2}{z_1}\right)
f_{i,j}\left(\frac{z_2}{z_1}\right)
T_{i+1}(x^{-1}z_1)T_j(z_2)-
f_{j,1}\left(\frac{x^{-i-1}z_1}{z_2}\right)f_{j,i}\left(\frac{z_1}{z_2}\right)
T_j(z_2)T_{i+1}(x^{-1}z_1)\notag\\
-&~c(r,x)\delta\left(\frac{x^{i-j}z_2}{z_1}\right)
f_{j+1,i}(x^{i-j+1})T_{j+1}(x^{j-i-1}z_1)T_i(z_1)\notag\\
+&~c(r,x)\delta\left(\frac{x^{i+j+2}z_2}{z_1}\right)
\prod_{l=1}^i \Delta_1(x^{2l+1})T_{i+j+1}(x^{i+1}z_2).
\notag
\end{align}
Using fusion relation (\ref{eqn:fusion1}) and 
$f_{j+1,i}(x^{i-j+1})T_{j+1}(x^{j-i-1}z_1)T_i(z_1)=f_{i,j+1}(x^{j-i-1})T_i(z_1)T_{j+1}(x^{j-i-1}z_1)$ in (\ref{induction4})
gives
\begin{align}
&
f_{i+1,j}\left(\frac{xz_2}{z_1}\right)T_{i+1}(x^{-1}z_1)T_j(z_2)-f_{j,i+1}\left(\frac{x^{-1}z_1}{z_2}\right)T_j(z_2)T_{i+1}(x^{-1}z_1)
\notag
\\
-&~c(r,x)\delta\left(\frac{x^{i-j}z_2}{z_1}\right)f_{i,j+1}(x^{-i+j+1})T_i(z_1)T_{j+1}(x^{-1}z_2)
\notag\\
+&~c(r,x)\delta\left(\frac{x^{i+j+2}z_2}{z_1}\right)\prod_{l=1}^{i}\Delta_1(x^{2l+1})T_{i+j+1}(x^{i+1}z_2).
\label{induction1}
\end{align}
Multiplying ${\rm RHS}_{i,j}(i)$ by $f_{1,i}\left(z_1/z_3\right)
f_{1,j}\left(z_2/z_3\right)T_1(z_3)$ from the left and using fusion relation (\ref{eqn:fusion2}) gives
\begin{align}
&c(r,x)\prod_{l=1}^{i-1}\Delta_1(x^{2l+1})
\left(\delta\left(\frac{x^{-i-j}z_2}{z_1}\right)f_{1,i+1}\left(\frac{x^jz_1}{z_3}\right)\Delta_1\left(\frac{x^i z_1}{z_3}\right)T_1(z_3)T_{i+j}(x^jz_1)\right.
\notag\\
-&~\left.
\delta\left(\frac{x^{i+j}z_2}{z_1}\right)f_{1,i+1}\left(\frac{x^{-j}z_1}{z_3}\right)\Delta_1\left(\frac{x^{-i} z_1}{z_3}\right)T_1(z_3)T_{i+j}(x^{-j}z_1)
\right).
\label{proof2}
\end{align}
Taking the limit $z_3\to x^{-i-1}z_1$ of (\ref{proof2}) multiplied by $c(r,x)^{-1}\left(1-x^{-i-1}z_1/z_3\right)$ 
and using fusion relation (\ref{eqn:fusion6}) along with the relation $\lim_{z_3 \to x^{-i-1}z_1}\left(1-x^{-i-1}z_1/z_3\right)
\Delta_1\left(x^{-i}z_1/z_3\right)=c(r,x)$ gives
\begin{align}
&
c(r,x)\delta\left(\frac{x^{-i-j}z_2}{z_1}\right)\prod_{l=1}^i \Delta_1(x^{2l+1})T_{i+j+1}(x^{-i-1}z_2)\notag\\
-~&
c(r,x)\delta\left(\frac{x^{i+j}z_2}{z_1}\right)
\prod_{l=1}^{i-1} \Delta_1(x^{2l+1})
f_{1,i+j}(x^{i-j+1})T_1(x^{-i-1}z_1)
T_{i+j}(x^{i}z_2).\label{induction2}
\end{align}
Multiplying ${\rm RHS}_{i,j}(k)$ $(1\leq k \leq i-1)$ by $f_{1,i}\left(z_1/z_3\right)f_{1,j}\left(z_2/z_3\right)T_1(z_3)$ 
from the left and using fusion relation (\ref{eqn:fusion3}) along with $f_{i-k,j+k}(x^{j-i})T_{i-k}(x^kz_1)T_{j+k}(x^{j-i+k}z_1)=f_{j+k,i-k}(x^{i-j})T_{j+k}(x^{j-i+k}z_1)T_{i-k}(x^kz_1)$ in (\ref{induction4})
gives
\begin{align}
&
c(r,x)\prod_{l=1}^{k-1}\Delta_1(x^{2l+1})
\label{proof3}
\\
\times&\left(
\delta\left(\frac{x^{-j+i-2k}z_2}{z_1}\right)
f_{1,i-k}\left(\frac{x^kz_1}{z_3}\right)f_{j+k,i-k}(x^{i-j})f_{1,j+k}\left(\frac{x^{-i+j+k}z_1}{z_3}\right)T_1(z_3)T_{j+k}(x^{j-i+k}z_1)T_{i-k}(x^kz_1)
\right.\notag\\
&-\left.
\delta\left(\frac{x^{j-i+2k}z_2}{z_1}\right)
f_{1,i-k}\left(\frac{x^{-k}z_1}{z_3}\right)f_{i-k,j+k}(x^{i-j})
f_{1,j+k}\left(\frac{x^{i-j-k}z_1}{z_3}\right)T_1(z_3)T_{i-k}(x^{-k}z_1)T_{j+k}(x^kz_2)
\right).\nonumber
\end{align}
Taking the limit $z_3\to x^{-i-1}z_1$ of (\ref{proof3}) multiplied by $c(r,x)^{-1}\left(1-x^{-i-1}z_1/z_3\right)$ 
and using fusion relations (\ref{eqn:fusion1}) and (\ref{eqn:fusion6}) along with
$$
f_{i-k+1,j+k}(x^{i-j+1})T_{i-k+1}(x^{-k-1}z_1)T_{j+k}(x^{-j+i-k}z_1)
=f_{j+k,i-k+1}(x^{j-i-1})T_{j+k}(x^{-j+i-k}z_1)T_{i-k+1}(x^{-k-1}z_1)
$$ in (\ref{induction4})
gives
\begin{align}
&
c(r,x)\prod_{l=1}^k \Delta_1(x^{2l+1})
\delta\left(\frac{x^{-j+i-2k}z_2}{z_1}\right)f_{j+k-1,i-k}(x^{i-j+1})
T_{i-k}(x^kz_1)T_{j+k+1}(x^{-k-1}z_2)\nonumber
\\
-~&
c(r,x)
\prod_{l=1}^{k-1}\Delta_1(x^{2l+1})
\delta\left(\frac{x^{j-i+2k}z_2}{z_1}\right)f_{i-k+1,j+k}(x^{i-j+1})T_{i-k+1}(x^{-k-1}z_1)T_{j+k}(x^kz_2).
\label{induction3}
\end{align}
Summing (\ref{induction1}), (\ref{induction2}), and (\ref{induction3}) for $1\leq k \leq i-1$ and shifting the variable $z_1 \mapsto xz_1$ gives ${\rm LHS}_{i+1,j}=\sum_{k=1}^{i+1}{\rm RHS}_{i+1,j}(k)$.
By induction on $i$, we have shown quadratic relation (\ref{CaseII-III:quadratic}).
~~$\Box$

\subsection{Classical limit}
\label{Section4}

The deformed $W$-algebra ${\cal W}_{q, t}({\mathfrak g})$ 
yields a $q$-Poisson $W$-algebra 
in the classical limit.
We set parameters $q=x^{2r}$ and $\beta=(r-1)/r$.
We obtain a $q$-Poisson $W$-algebra \cite{Frenkel-Reshetikhin1,Frenkel-Reshetikhin-Semenov,Semenov-Sevostyanov} 
in the classical limit $\beta \to 0$ with $q$ fixed. 
The defining relations of the deformed $W$-superalgebra ${\cal W}_{q, t}(\mathfrak{sl}(2|1))$ are given by
\begin{align}
[T_i[m],T_j[n]]
&=-\sum_{l=1}^\infty f_{i,j}^l
\left(T_i[m-l]T_j[n+l]-T_j[n-l]T_i[m+l]\right)\notag
\\
&+c(r,x)\sum_{k=1}^i \prod_{l=1}^{k-1}\Delta_1(x^{2l+1})\notag\\
&
\times \sum_{l \in {\mathbf Z}}\left(
f_{i-k,j+k}(x^{j-i})x^{(j-i)l+k(m-n)+4kl}T_{i-k}[m-l]T_{j+k}[n+l]
\right.\notag
\\
&
\left.-
f_{i-k,j+k}(x^{-j+i})x^{(i-j)l+k(n-m)-4kl}T_{i-k}[m-l]T_{j+k}[n+l]
\right).\notag
\end{align}
We define $f_{i,j}^l$ by $f_{i,j}(z)=\sum_{l=0}^\infty f_{i,j}^l z^l$, where the structure functions $f_{i,j}(z)$ are given in (\ref{CaseI-II-III:fij}).
We define the $q$-Poisson bracket $\{,\}$ by taking the classical limit $\beta \to 0$ with $q$ fixed as
\begin{eqnarray}
\{T_i^{{\rm PB}}[m], T_j^{{\rm PB}}[n]\}=-\lim_{\beta \to 0}
\frac{1}{\beta \log q}[T_i[m],T_j[n]].\nonumber
\end{eqnarray}
Here, we set $T_i^{PB}[m]$ as
\begin{eqnarray}
T_i(z)=\sum_{m \in {\mathbf Z}}T_i[m]z^{-m} \longrightarrow T_i^{PB}(z)=\sum_{m \in {\mathbf Z}}T_i^{PB}[m]z^{-m}
~~(\beta \to 0,~q~\rm{fixed}).\nonumber
\end{eqnarray}
The $\beta$-expansions of the structure functions are given as 
\begin{eqnarray}
&&f_{i,j}(z)=1+\beta \log q \sum_{m=1}^\infty
\frac{\left[\frac{1}{2} {\rm Min}(i,j)m \right]_q \left[\left(\frac{1}{2}{\rm Max}(i,j)-1\right)m\right]_q}{[m]_q} (q-q^{-1})+O(\beta^2)~~(i, j \geq 1),
\nonumber\\
&&c(r,x)=-\beta \log q+O(\beta^2).
\nonumber
\end{eqnarray}
We obtain the following Proposition.
\begin{prop}
For the $q$-Poisson $W$-superalgebra for $\mathfrak{sl}(2|1)$
the generating function
$T_i^{PB}(z)$ satisfies
\begin{align}
\{T_i^{PB}(z_1),T_j^{PB}(z_2)\}
&=(q-q^{-1})C_{i,j}\left(\frac{z_2}{z_1}\right)
T_i^{PB}(z_1)T_j^{PB}(z_2)
\notag
\\
&+\sum_{k=1}^i \delta\left(\frac{q^{\frac{-j+i}{2}-k}z_2}{z_1}\right)T_{i-k}^{PB}(q^{\frac{k}{2}}z_1)T_{j+k}^{PB}(q^{-\frac{k}{2}}z_2)
\notag
\\
&-\sum_{k=1}^i \delta\left(\frac{q^{\frac{j-i}{2}+k}z_2}{z_1}\right)T_{i-k}^{PB}(q^{-\frac{k}{2}}z_1)T_{j+k}^{PB}(q^{\frac{k}{2}}z_2)~~(1\leq i \leq j).
\notag
\end{align}
Here $T_0^{PB}(z)$ in the right hand side is understood as $T_0^{PB}(z)=1$.
Here we set the structure functions $C_{i,j}(z)$ $(i,j\geq 1)$ as
\begin{eqnarray}
C_{i,j}(z)=\sum_{m \in {\mathbf Z}}
\frac{\left[\frac{1}{2} {\rm Min}(i,j)m \right]_q \left[\left(\frac{1}{2}{\rm Max}(i,j)-1\right)m\right]_q}{[m]_q} z^m~~(i,j \geq 1).
\nonumber
\end{eqnarray}
\end{prop}

The structure functions satisfy
$C_{i,2}(z)=C_{2,i}(z)=0$ $(i=1,2)$.

\section{Conclusion and Discussion}
\label{Section4}

We revisited Ding-Feigin's construction of the deformed $W$-algebras
${\cal W}_{q,t}(\mathfrak{sl}(3))$ and ${\cal W}_{q,t}(\mathfrak{sl}(2|1))$.
Using Ding-Feigin realization we introduced the higher $W$-currents $T_i(z)$
for the deformed $W$-superalgebra ${\cal W}_{q, t}(\mathfrak{sl}(2|1))$.
We obtained a set of quadratic relations of $T_i(z)$,
which is independent of the choice of Dynkin diagrams for the superalgebra $\mathfrak{sl}(2|1)$,
though the screening currents depend on it.
There is an infinite number of quadratic relations
for an infinite number of the
currents $T_i(z)$
for ${\cal W}_{q,t}(\mathfrak{sl}(2|1))$.
The higher $W$-currents $T_i(z)$ can also be constructed by repeating fusion 
from $T_1(z)$ as
\begin{eqnarray}
T_i(z)=\frac{1}{c(r,x)}
{\displaystyle
\lim_{w \to x^{-1}z}
\left(1-\frac{x^{-1}z}{w}\right) 
f_{i-1,1}\left(\frac{x^{i-1}z}{w}\right)}
T_{i-1}(w)T_1(x^{i-1}z)~~~(i \geq 2),
\label{fusion:Ti(z)-T1(z)}
\end{eqnarray}
which is a special case of lemma \ref{lem:III-4}.
For ${\cal W}_{q,t}(\mathfrak{sl}(3))$, as for ${\cal W}_{q,t}(\mathfrak{sl}(2|1))$, 
the higher $W$-currents $T_i(z)$
can be formally defined by (\ref{fusion:Ti(z)-T1(z)}).
However, with the truncations
\begin{eqnarray}
{\displaystyle \lim_{z_1 \to x^{-2}z_2}\left(
1-\frac{x^{-2}z_2}{z_1}\right)}f_{1,1}\left(\frac{z_2}{z_1}\right)\Lambda_k(z_1)\Lambda_k(z_2)=0~~
(1\leq k \leq 3),~~
:\Lambda_1(x^{-2}z)\Lambda_2(z)\Lambda_3(x^2z):=1,
\nonumber
\end{eqnarray}
we have $T_3(z)=1$ and $T_i(z)=0$ $(i \geq 4)$. 
We get only $T_1(z)$ and $T_2(z)$ from (\ref{fusion:Ti(z)-T1(z)}) for ${\cal W}_{q,t}(\mathfrak{sl}(3))$.\\

It seems to be possible to extend Ding-Feigin construction to the case of many fermions which will give
a higher rank generalization ${\cal W}_{q, t}(\mathfrak{sl}(M|N))$, and obtain its quadratic relations. 
We expect to report it in the near future. 
It is still an open problem to find quadratic relations of the deformed $W$-algebra ${\cal W}_{q,t}({\mathfrak g})$,
except for ${\mathfrak g}=\mathfrak{sl}(N)$, $A_2^{(2)}$, and $\mathfrak{sl}(2|1)$.
It seems to be possible to extend Ding-Feigin construction to other superalgebras and obtain
their quadratic relations. 

~\\
{\bf ACKNOWLEDGMENTS}

This paper is dedicated to Professor Boris Feigin on the occasion of his 65th birthday.
The author would like to thank Professor Michio Jimbo for carefully reading the manuscript many times and for giving lots of useful advice.
The author would like to thank Professor Kenji Iohara and Professor Ryusuke Endo for discussions.
This work is supported by the Grant-in-Aid for Scientific Research {\bf C} (26400105) and {\bf C} (19K03509) from the Japan Society for the Promotion of Science.

\begin{appendix}
~\\
\begin{Large}{\bf Appendix}\end{Large}
\renewcommand{\theequation}{\Alph{section} \arabic{equation}}
\setcounter{equation}{0}
\section{Fusion relations}
\label{Appendix:fusion}
In this appendix we summarize the fusion relations of $\Lambda_i(z)$ for Case 2
which are obtained from Proposition \ref{prop:III-3}.
The relations for Case 3 are given in the same way.
We use the abbreviations
\begin{eqnarray}
&&
\Lambda_3^{(i)}(z)=:\prod_{l=1}^i \Lambda_3(x^{-i-1+2l}z):~~~(i \geq 1),
\nonumber
\\
&&
\Lambda_{k,3}^{(i)}(z)=
\left\{\begin{array}{cc}
:\Lambda_k(x^{-i+1}z) {\displaystyle \prod_{l=1}^{i-1}}\Lambda_3(x^{-i+1+2l}z):& (i \geq 2),\\
\Lambda_k(z)& (i=1)
\end{array}\right.~~~(1\leq k \leq 2),
\label{abbreviation:lambda}
\end{eqnarray}
\begin{eqnarray}
&&
\Lambda_{1,2,3}^{(i)}(z)=
\left\{\begin{array}{cc}
{\displaystyle :\Lambda_1(x^{-i+1}z)\Lambda_2(x^{-i+3}z)\prod_{l=1}^{i-2}\Lambda_3(x^{-i+3+2l}z)}:&(i \geq 3),\\
:\Lambda_1(x^{-1}z)\Lambda_2(xz):&(i=2),\\
0&(i=1).
\end{array}
\right.\nonumber
\end{eqnarray}
We set
\begin{eqnarray}
F_{i,j}^{(\pm)}(z)=(1-x^{\pm (i+j)}z)f_{i,j}(z).\nonumber
\end{eqnarray}
%%%%%%%%%%%%%%%%%%%%%%%%%%%%%%%%%%%%%
$\bullet$~For $k=1,2$, we have
\begin{eqnarray}
\lim_{z_1 \to x^{\pm (i+j)}z_2}
F_{i,j}^{(\pm)}\left(\frac{z_2}{z_1}\right)\Lambda_3^{(i)}(z_1)\Lambda_{3}^{(j)}(z_2)=
\mp \frac{c(r,x)d_{i+j}(r,x)}{d_{i}(r,x)d_{j}(r,x)}\prod_{l=1}^{{\rm Min}(i,j)-1}\Delta_1(x^{2l+1})\Lambda_{3}^{(i+j)}(x^{\pm i}z_2)\nonumber\\
(i,j \geq 1),
\label{appendixA:1}
\\
\lim_{z_1 \to x^{i+j}z_2}
F_{i,j}^{(+)}\left(\frac{z_2}{z_1}\right)\Lambda_3^{(i)}(z_1)\Lambda_{1,2,3}^{(j)}(z_2)
=
-\frac{c(r,x)d_{i+j-2}(r,x)}{d_{i}(r,x)d_{j-2}(r,x)}\prod_{l=1}^{{\rm Min}(i,j)-1}\Delta_1(x^{2l+1})\Lambda_{1,2,3}^{(i+j)}(x^{i}z_2)\nonumber\\
(i \geq 1, j \geq 2),
\label{appendixA:2}
\\
\lim_{z_1 \to x^{i+j}z_2}
F_{i,j}^{(+)}\left(\frac{z_2}{z_1}\right)\Lambda_3^{(i)}(z_1)\Lambda_{k,3}^{(j)}(z_2)
=
-\frac{c(r,x)d_{i+j-1}(r,x)}{d_{i}(r,x)d_{j-1}(r,x)}\prod_{l=1}^{{\rm Min}(i,j)-1}\Delta_1(x^{2l+1})\Lambda_{k,3}^{(i+j)}(x^{i}z_2)\nonumber
\\
(i,j \geq 1),\label{appendixA:3}
\\
\lim_{z_1 \to x^{-(i+j)}z_2}
F_{i,j}^{(-)}\left(\frac{z_2}{z_1}\right)\Lambda_{1,2,3}^{(i)}(z_1)\Lambda_{3}^{(j)}(z_2)=
\frac{c(r,x)d_{i+j-2}(r,x)}{d_{i-2}(r,x)d_{j}(r,x)}\prod_{l=1}^{{\rm Min}(i,j)-1}\Delta_1(x^{2l+1})\Lambda_{1,2,3}^{(i+j)}(x^{-i}z_2)\nonumber
\\
(i \geq 2, j \geq 1),
\label{appendixA:4}
\\
\lim_{z_1 \to x^{-(i+j)}z_2}
F_{i,j}^{(-)}\left(\frac{z_2}{z_1}\right)
\Lambda_{k,3}^{(i)}(z_1)\Lambda_{3}^{(j)}(z_2)
=\frac{c(r,x)d_{i+j-1}(r,x)}{d_{i-1}(r,x)d_{j}(r,x)}\prod_{l=1}^{{\rm Min}(i,j)-1}\Delta_1(x^{2l+1})\Lambda_{k,3}^{(i+j)}(x^{-i}z_2)\nonumber\\
(i,j \geq 1). \label{appendixA:5}
\end{eqnarray}
$\bullet$~As exceptional formulae
of $\Lambda_3^{(i)}(z)$, $\Lambda_{k,3}^{(i)}(z)$, and $\Lambda_{1,2,3}^{(i)}(z)$
for small $i$, we have
\begin{eqnarray}
\lim_{z_1 \to x^{(j+1)}z_2}
F_{j,1}^{(+)}\left(\frac{z_2}{z_1}\right)
\Lambda_{2,3}^{(j)}(z_1)\Lambda_{1}(z_2)=
-c(r,x)\Lambda_{1,2,3}^{(j+1)}(x^{j}z_2)~~~(j \geq 2),
\label{appendixA:6}
\\
\lim_{z_1 \to x^{-(j+1)}z_2}
F_{1,j}^{(-)}\left(\frac{z_2}{z_1}\right)
\Lambda_{1}(z_1)\Lambda_{2,3}^{(j)}(z_2)=
c(r,x)\Lambda_{1,2,3}^{(j+1)}(x^{-1}z_2)~~~(j \geq 2),
\label{appendixA:7}
\end{eqnarray}
\begin{eqnarray}
&&
\lim_{z_1 \to x^{\pm 2}z_2}F_{1,1}^{(\pm)}\left(\frac{z_2}{z_1}\right)\Lambda_k(z_1)\Lambda_k(z_2)=0~~~(k=1,2).
\label{appendixA:8}
\end{eqnarray}
For $1\leq k<l\leq 3$, we have
\begin{eqnarray}
&&
\lim_{z_1 \to x^{2}z_2}
F_{1,1}^{(+)}\left(\frac{z_2}{z_1}\right)\Lambda_k(w_1)\Lambda_l(w_2)=0,
~~
\lim_{z_1 \to x^{-2}z_2}
F_{1,1}^{(-)}\left(\frac{z_2}{z_1}\right)\Lambda_l(z_1)\Lambda_k(z_2)=0,
\label{appendixA:10}
\\
&&
\lim_{z_1 \to x^{2}z_2}
F_{1,1}^{(+)}\left(\frac{z_2}{z_1}\right)\Lambda_l(z_1)\Lambda_k(z_2)
=-c(r,x):\Lambda_k(z_2)\Lambda_l(x^2z_2):,
\label{appendixA:11}
\\
&&
\lim_{z_1 \to x^{-2}z_2}
F_{1,1}^{(-)}\left(\frac{z_2}{z_1}\right)\Lambda_k(z_1)\Lambda_l(z_2)
=c(r,x):\Lambda_k(x^{-2}z_2)\Lambda_l(z_2):.\label{appendixA:12}
\end{eqnarray}

The remaining fusions vanish.
\section{Exchange relations}
\label{Appendix:exchange}
\setcounter{equation}{0}

In this appendix we give the exchange relations of $\Lambda_i(z)$ for Case 2
which are obtained from Proposition \ref{prop:III-3}, (\ref{eqn:fusion5}), and (\ref{eqn:fusion4}).
The relations for Case 3 are obtained in the same way.
We use the abbreviations (\ref{abbreviation:lambda}) and set
\begin{eqnarray}
{f_{i,j}}\left(
{\cal A}(z_1), {\cal B}(z_2)\right)_{1 \leftrightarrow 2}
=f_{i,j}\left(\frac{z_2}{z_1}\right){\cal A}(z_1){\cal B}(z_2)-
f_{j,i}\left(\frac{z_1}{z_2}\right){\cal B}(z_2){\cal A}(z_1).
\nonumber
\end{eqnarray}

For $i \geq 1$ and $k=1,2$, we have
\begin{eqnarray}
&&
f_{1,i}\left(\Lambda_k(z_1), \Lambda_{1,2,3}^{(i)}(z_2)\right)_{1\leftrightarrow 2}=0,~~~
f_{1,i}\left(\Lambda_k(z_1), \Lambda_{k,3}^{(i)}(z_2)\right)_{1\leftrightarrow 2}=0,
\label{appendixB:1}
\end{eqnarray}
\begin{eqnarray}
&&
f_{1,i}\left(\Lambda_1(z_1), \Lambda_{2,3}^{(i)}(z_2)\right)_{1 \leftrightarrow 2}
=c(r,x)
\left(\delta\left(\frac{x^{-i-1}z_2}{z_1}\right)-\delta\left(\frac{x^{-i+1}z_2}{z_1}\right)\right)
:\Lambda_1(z_1)\Lambda_{2,3}^{(i)}(z_2):,
\label{appendixB:2}
\end{eqnarray}
\begin{eqnarray}
f_{1,i}\left(\Lambda_k(z_1),\Lambda_{3}^{(i)}(z_2)\right)_{1\leftrightarrow 2}
=c(r,x)
\left(\delta\left(\frac{x^{-i-1}z_2}{z_1}\right)-\delta\left(\frac{x^{-i+1}z_2}{z_1}\right)\right)
:\Lambda_k(z_1)\Lambda_{3}^{(i)}(z_2):,
\label{appendixB:3}
\end{eqnarray}
\begin{eqnarray}
f_{1,i}\left(\Lambda_3(z_1), \Lambda_{3}^{(i)}(z_2)\right)_{1\leftrightarrow 2}
=\frac{c(r,x)d_{i+1}(r,x)}{d_1(r,x)d_i(r,x)}
\left(\delta\left(\frac{x^{-i-1}z_2}{z_1}\right)-\delta\left(\frac{x^{i+1}z_2}{z_1}\right)\right)
:\Lambda_3(z_1)\Lambda_{3}^{(i)}(z_2):,
\label{appendixB:4}
\end{eqnarray}
\begin{eqnarray}
f_{1,i}\left(\Lambda_2(z_1), \Lambda_{1,3}^{(i)}(z_2)\right)_{1\leftrightarrow 2}=
c(r,x)
\left(\delta\left(\frac{x^{-i+1}z_2}{z_1}\right)-\delta\left(\frac{x^{-i+3}z_2}{z_1}\right)\right)
:\Lambda_2(z_1)\Lambda_{1,3}^{(i)}(z_2):,
\label{appendixB:5}
\end{eqnarray}
\begin{eqnarray}
f_{1,i}\left(\Lambda_3(z_1),\Lambda_{k,3}^{(i)}(z_2)\right)_{1\leftrightarrow 2}
=\frac{c(r,x)d_{i}(r,x)}{d_1(r,x)d_{i-1}(r,x)}
\left(\delta\left(\frac{x^{-i+1}z_2}{z_1}\right)-\delta\left(\frac{x^{i+1}z_2}{z_1}\right)\right)
:\Lambda_3(z_1)\Lambda_{k,3}^{(i)}(z_2):.
\label{appendixB:6}
\end{eqnarray}

For $i \geq 2$, we have
\begin{eqnarray}
f_{1,i}\left(\Lambda_3(z_1),\Lambda_{1,2,3}^{(i)}(z_2)\right)_{1\leftrightarrow 2}=
\frac{c(r,x)d_{i-1}(r,x)}{d_1(r,x)d_{i-2}(r,x)}
\left(\delta\left(\frac{x^{-i+3}z_2}{z_1}\right)-\delta\left(\frac{x^{i+1}z_2}{z_1}\right)\right)
:\Lambda_3(z_1)\Lambda_{1,2,3}^{(i)}(z_2):.
\label{appendixB:7}
\end{eqnarray}

\end{appendix}
\end{document}